
\input amstex

\magnification=1200
\loadmsam
\loadmsbm
\loadeufm
\loadeusm
\UseAMSsymbols

\hsize=6.0truein
\hoffset=0.15truein
\vsize=9truein
\voffset=-0.2truein

\def\leftitem#1{\item{\hbox to\parindent{\enspace#1\hfill}}}

\def\boxit#1#2{\hbox{\vrule
	\vtop{%
	\vbox{\hrule\kern#1%
	\hbox{\kern#1#2\kern#1}}%
	\kern#1\hrule}%
	\vrule}}

\def\leaderfill{\leaders\hbox to 1em{\hss.\hss}\hfill}

\parskip=\medskipamount
\document

\input epsf

\centerline{\bf Simple Loops on Surfaces and Their Intersection Numbers}
\bigskip
\centerline{Feng Luo}

\centerline{\it Dept. of Math., Rutgers University,  New Brunswick, NJ 08903 \rm}
\centerline{e-mail: fluo\@math.rutgers.edu}

{\bf Abstract.} 
Given a compact orientable surface $\Sigma$, let  $\Cal S(\Sigma)$ be 
the set of isotopy classes
of essential simple loops on $\Sigma$. We  determine a complete set of relations
for a function from $\Cal S(\Sigma)$ to $\bold Z$  to be 
a geometric intersection
number function. As a consequence, we obtain explicit equations in
$\bold R^{\Cal S(\Sigma)}$ and $P (\bold R^{\Cal S(\Sigma)})$
defining Thurston's space of measured laminations and Thurston's 
compactification of the Teichm\"uller space. These equations are not
only piecewise integral linear but also semi-real algebraic.

Table of contents:

\S 1. Introduction

\S 2.  A Multiplicative Structure on Curve Systems

\S 3. The One-holed Torus

\S 4. The Four-holed Sphere

\S 5. A Reduction Proposition

\S 6. The Two-holed Tours and the  Five-holed Sphere

\S 7. Proofs of Theorem 1 and the Corollary

\S 8. Proofs of Results in Section 2 and Some Questions

References

\bigskip

\S 1. Introduction

Given a compact orientable surface $\Sigma$ =$\Sigma_{g,r}$ of genus $g$ with $r$
boundary components, let $\Cal S$ = $\Cal S(\Sigma)$ be the set of isotopy classes
of essential simple loops on $\Sigma$. A function $f: \Cal S(\Sigma) \to \bold R$
is called a \it geometric intersection number function, \rm or simply
\it geometric function \rm if there is a  measured lamination $m$ on $\Sigma$ so that
$f(\alpha)$ is the measure of $\alpha$ in $m$. Geometric functions were
introduced and studied by W. Thurston in his work  on the  classification
of surface homeomorphisms and the compactification  of the Teichm\"uller 
spaces ([FLP], [Th]). The space of all geometric functions under the 
pointwise convergence topology is homeomorphic to Thurston's space of
measured laminations $\Cal ML(\Sigma)$.  Thurston showed that  $\Cal ML(\Sigma)$ is homeomorphic to
a Euclidean space and $\Cal ML(\Sigma)$ has a piecewise integral linear structure invariant
under the action of the mapping class group. The projectivization of 
$\Cal ML(\Sigma)$ is Thurston's boundary of the Teichm\"uller space. The object of the paper
is to characterize all geometric functions on $\Cal S(\Sigma)$. As a consequence, both
$\Cal ML(\Sigma)$ and its projectivization are reconstructed explicitly in 
terms of an intrinsic $(\bold QP^1, PSL(2, \bold Z))$ structure on $\Cal S(\Sigma)$.

{\bf Theorem 1.} \it Suppose $\Sigma$ is a compact orientable surface of  negative
Euler number. Then a function $f$ on $\Cal S(\Sigma)$ is geometric if and only if for
each incompressible subsurface $\Sigma' \cong \Sigma_{1,1}$ or $\Sigma_{0,4}$,
the restriction $f|_{\Cal S(\Sigma')}$ is geometric. Furthermore, geometric
functions on $\Cal S(\Sigma_{1,1})$ and $\Cal S(\Sigma_{0,4})$ are 
characterized  by two homogeneous equations in the $(\bold QP^1, PSL(2, \bold Z))$ structure on $\Cal S(\Sigma)$. \rm

Recall that a  subsurface $\Sigma' \subset$ $\Sigma$ is \it incompressible \rm if each essential
loop in $\Sigma'$ is still essential in  $\Sigma$. It is well known that if each boundary
component of $\Sigma'$ is essential in $\Sigma$, then $\Sigma'$ is essential.

Geometric functions and measures laminations haven been studied from many
different points of views. Especially, they are identified with height
functions and horizontal foliations associated to  holomorphic quadratic
forms on $\Sigma$ ([Ga], [HM], [Ker1]). They are also related to the translation
length functions of group action on $\bold R$-trees and the valuation theory
([Bu], [CM], [MS], [Par]). 
In [Bo1], measured laminations and hyperbolic metrics are
considered as special cases of currents. As a consequence, 
Thurston's compactification is derived from a natural setting. 

Our approach is  combinatorial and is based on the notion of curve systems
([De], [FLP], [Hat], [PH],  [Th]). Recall that 
a \it curve system \rm is a finite
disjoint union of essential proper arcs and essential non-boundary parallel
simple loops on the surface. Let $\Cal CS(\Sigma)$ be the set of isotopy 
classes of curve systems on $\Sigma$. The space  $\Cal CS(\Sigma)$ was 
introduced by Dehn and rediscovered independently by Thurston. Dehn  called 
the space  the \it arithmetic field \rm of the topological surface.
Given two classes $\alpha$, $\beta$ in $\Cal CS(\Sigma) \cup \Cal S(\Sigma)$, 
their \it geometric intersection number \rm
$I(\alpha, \beta)$ is defined to be min$\{|a \cap b|: a \in \alpha, b \in
\beta\}$. The essential part of the paper is to characterize those 
geometric functions $f$ so that $f(\alpha) = I(\alpha, \beta) 
(=I_{\beta}(\alpha))$ for some fixed $\beta \in$ $\Cal CS(\Sigma)$.  

\midspace{0.1cm}
\centerline{\epsfbox{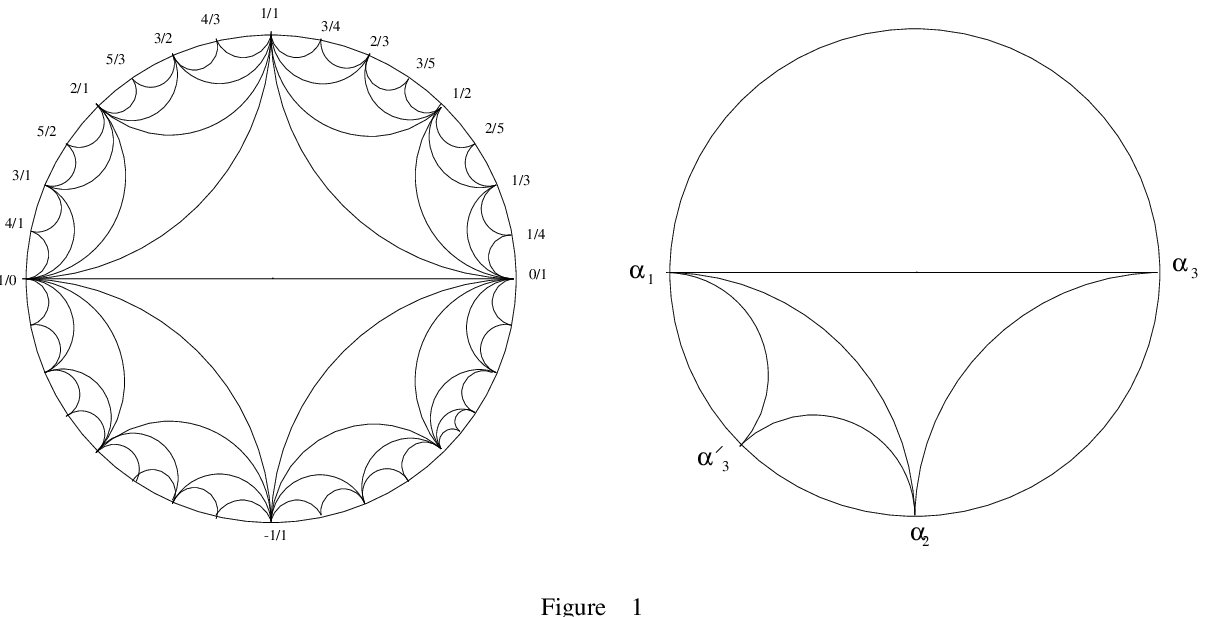}}
\midspace{0.1cm}

Given a surface $\Sigma$, let $\Cal S'(\Sigma)$ = 
$\Cal CS(\Sigma)$$\cap$$\Cal S(\Sigma)$ be the set of isotopy classes of
essential, non-boundary parallel simple loops in $\Sigma$. For surfaces $\Sigma$ = $\Sigma_{1,0}$, $\Sigma_{1,1}$ and $\Sigma_{0,4}$, it is well 
known that there exists
a bijection $\pi: \Cal S'(\Sigma) \to \bold QP^1 (= \hat \bold Q) $ 
so that $p'q - pq' =
\pm 1$ if and only if $I(\pi^{-1}(p/q), \pi^{-1}(p'/q')) =1$ 
(for $\Sigma_{1,0}$, $\Sigma_{1,1}$)
and 2 (for $\Sigma_{0,4}$). See figure 1.
We say that three distinct classes $\alpha$, $\beta$, $\gamma$ in 
$\Cal S'(\Sigma)$ form an \it
ideal triangle \rm if they correspond to the vertices of an ideal
triangle in the modular relation under the map $\pi$.

{\bf Theorem 2.} \it (a) For surface $\Sigma_{1,1}$, a function $f: \Cal S$
$\to \bold Z_{\geq 0}$  is a geometric
function $I_{\delta}$ with $\delta \in $$\Cal CS(\Sigma)$ if and only if
the following hold. 
$$
f(\alpha_1) + f(\alpha_2) + f(\alpha_3) = \max_{i=1,2,3}(2 f(\alpha_i),
f([\partial \Sigma_{1,1}])) \tag 1
$$ where $(\alpha_1, \alpha_2, \alpha_3)$
 is an ideal triangle, and
$$f(\alpha_3) + f(\alpha_3') = \max(2f(\alpha_1), 2f(\alpha_2),
f([\partial \Sigma_{1,1}])) \tag 2 $$ 
where   $(\alpha_1, \alpha_2, \alpha_3)$ and
 $(\alpha_1, \alpha_2, \alpha_3')$  are two distinct ideal triangles.
 $$f([\partial \Sigma_{1,1}]) \in 2\bold Z. \tag 3$$

(b) For surface $\Sigma_{0,4}$ with $\partial \Sigma_{0,4} = b_1 \cup b_2 \cup
b_3 \cup b_4$, a function $f: \Cal S \to \bold Z_{\geq 0}$ is a geometric
function $I_{\delta}$ for some $\delta \in$ $\Cal CS(\Sigma)$ if and only if
for each ideal triangle $(\alpha_1, \alpha_2, \alpha_3)$ so that
$(\alpha_i, b_s, b_r)$ bounds a $\Sigma_{0,3}$ in $\Sigma_{0,4}$ the following
hold.
 $$\Sigma_{i=1}^3 f(\alpha_i) = \max_{1 \leq i \leq 3;1 \leq s \leq 4}
(2f(\alpha_i), 2f(b_s), \sum_{s=1}^4 f(b_s), f(\alpha_i)+ f(b_s) + f(b_r))
\tag 4$$
$$f(\alpha_3) + f(\alpha_3') =\max_{1 \leq i \leq 2;1 \leq s \leq 4}
(2f(\alpha_i), 2f(b_s), \sum_{s=1}^4 f(b_s), f(\alpha_i)+ f(b_s) + f(b_r))
\tag 5$$
where   $(\alpha_1, \alpha_2, \alpha_3)$ and
$(\alpha_1, \alpha_2, \alpha_3')$  are two distinct ideal triangles,
$$f(\alpha_i) + f(b_s) + f(b_r) \in 2 \bold Z. \tag 6$$

(c) The characterization of geometric functions $f: \Cal S(\Sigma) \to \bold 
R_{\geq 0}$ for $\Sigma= \Sigma_{1,1}$ and $\Sigma_{0,4}$  is  given
by equations (1),(2) (for $\Sigma_{1,1}$) and (4), (5) (for $\Sigma_{0,4}$).\rm

Theorem 2 is motivated by the tours case.  In fact for the torus $\Sigma_{1,0}$,
a function on $\Cal S(\Sigma_{1,0})$ is geometric if and only if it satisfies the triangular equality
$f(\alpha_1) + f(\alpha_2) + f(\alpha_3) = \max_{i=1,2,3}(2 f(\alpha_i))$
and  $f(\alpha_3) + f(\alpha_3') =\max(2f(\alpha_1), 2f(\alpha_2))$.

The equations (1),(2),(4) and (5) in theorem 2 are obtained as the
degenerations of the trace identities for $SL(2, \bold R)$ matrices.
For instance, equations (1), (2) are the  degenerations of $tr(A)tr(B)tr(AB)
= tr^2(A) + tr^2(B) + tr^2(AB) - tr([A,B]) -2$ and
$tr(AB) tr(A^{-1}B) = tr^2(A) + tr^2(B) -tr([A,B]) -2$.
 
Several properties of the measured laminations spaces are reflected
in the equations (1),(2),(4), and (5). For instance, since
the equations are piecewise integral linear so that 
rational solutions are dense,  one obtains Thurston's result
that  the space $\Cal ML(\Sigma)$ has a piecewise integral linear structure 
and the rational multiples of the  curve systems is a  dense subset. 
On the other hand, the equations are also 
semi-real algebraic. 
Indeed, the space defined by
$\sum_{i=1}^k x_i = \max_{1 \leq j \leq l} (y_j)$ is  semi-real
algebraic since it is  equivalent to:
$\prod_{j=1}^l( \sum_{i=1}^k x_i - y_j) = 0$, and
$\sum_{i=1}^k x_i \geq y_j$, for all $j$. This seems to indicate that
the space $\Cal ML(\Sigma)$ has a semi-real algebraic structure. Given a
surface $\Sigma_{g,r}$, Thurston showed that there exists a finite set
$F$ consisting of $9g+4r-9$ elements in $\Cal S(\Sigma)$ so that the map $\tau_F:$ $\Cal ML(\Sigma)$
$\to \bold R_{\geq 0}^F$ sending $m$ to $I_m |_F$ is an embedding ([FLP]). As a
consequence of theorems 1,2, we have,

{\bf Corollary}. \it For surface $\Sigma_{g,r}$ of negative Euler number,
there is a finite set $F$ consisting of $9g+4r-9$ elements in $\Cal S(\Sigma)$ so that
the map $\tau_F$ is an embedding whose image is a polyhedron defined
by finitely many explicit integer coefficient polynomial equations and
inequalities. \rm

It is interesting to observe  that the approach taken in the
paper (also in [Lu1], [Lu3])  follows Grothendieck's philosophy of the
``Teichm\"uller tower" where the ``generators" are the surfaces $\Sigma_{1,1}$
and $\Sigma_{0,4}$ and the ``relations" are $\Sigma_{1,2}$ and $\Sigma_{0,5}$.
See [Sch] for more details.
From this point of view, it seems clear that  the
$(\bold QP^1, PSL(2, \bold  Z))$ modular structure
is fundamental to the topology and geometry of surfaces and the
modular structure plays a role of ``local coordinate" on the set
$\Cal S(\Sigma)$. Following this line,
we may ask the following two questions on the
related topics of mapping class groups and SL(2,$\bold C)$ representations.

Question 1. (A presentation of the mapping class group).
Suppose $\Sigma$ is a compact oriented surface. Let $Mod(\Sigma)$
be the mapping class group of $\Sigma$ consisting of  isotopy classes of
orientation preserving homeomorphisms which leaves each boundary
component invariant. Let $G$ be the group with $\Cal S(\Sigma)$
as the set of generators and the following as the set of relations:
$(R_1)$ $xy=yx$ if $I(x,y)=0$;
$(R_2)$  $x = 1$ if $x$ is a boundary component of $\Sigma$;
$(R_3)$ $xy=yz$ if ($x,y,z$) forms 
a positively oriented  ($x \to y \to z \to x$ is
the right hand order in $S^1$) ideal triangle
in $\Cal S(\Sigma')$ where $ \Sigma' \cong \Sigma_{1,1}$ is
incompressible in $\Sigma$;
$(R_4)$ $xyz = b_1b_2b_3b_4$ if ($x,y,z$) forms a positively oriented
ideal triangle in $\Cal S(\Sigma')$ where $\Sigma' \cong \Sigma_{0,4}$ is
incompressible in $\Sigma$ with $\partial \Sigma' = b_1 \cup
b_2 \cup b_3 \cup b_4$.
Is $G$ a presentation of $Mod(\Sigma)$?

Note that relation $(R_3)$ implies the Artin's relation ($xyx=yxy$) 
and $(R_4)$ is the lantern relation which was discovered
by Dehn ([De], p333) in 1938 and rediscovered independently by 
Johnson. See [Bi], [De], [Har], [HT], [Li], [Waj] for more details.

Question 2. (Characters of SL(2,$\bold C)$ representations)
A function $f : \Cal S(\Sigma) \to \bold C$ is the (restriction of) character of a 
representation of $\pi_1(\Sigma)$ into $SL(2, \bold C$) if and only
if $f|_{\Cal S(\Sigma')}$ is a character for each incompressible
subsurface $\Sigma' \cong \Sigma_{1,1}$ or $\Sigma_{0,4}$.

The description of characters for  the surfaces $\Sigma_{1,1}$ and
$\Sigma_{0,4}$ seems to be known.
See [CS], [Go], [GoM], [Ho], [Mag] and the references cited therein.  

The organization of the paper is as  follows. In \S2, we establish several
basic properties of the curve systems. In particular, a multiplicative
structure on $\Cal CS(\Sigma)$ is introduced. In \S3,\S4, we prove
theorem 2. The proof in \S4 is complicated due to the existence of eight
different ideal triangulations of the surface $\Sigma_{0,4}$.
In \S5, we prove a reduction result. This is one of the key steps in the
proof of theorem 1. It reduces the general case to two
surfaces: $\Sigma_{1,2}$ and $\Sigma_{0,5}$. In \S6, we prove 
theorem 1 for surfaces $\Sigma_{1,2}$
and $\Sigma_{0,5}$. The proofs of  theorem 1 and the corollary are in
\S7. The proof of the results in \S2 is in  \S8.

\it Acknowledgment. \rm I  would like to thank F. Bonahon, M. Freedman,
X.S. Lin,  and Y. Minsky for discussions. The work 
is supported in part by the NSF.

\S2.  A Multiplicative Structure on Curve Systems

We work in the piecewise linear category. Surfaces are oriented and connected
and  have negative Euler numbers unless specified otherwise.
A regular neighborhood of a submanifold $X$ is denoted by $N(X)$. Regular
neighborhoods are assumed to be small. The isotopy class of a curve system $c$
will be denoted by $[c]$. Suppose $f:$ $\Cal CS(\Sigma)$ $\to \bold R$ is a function and
$c$ is a curve system. We define $f(c)$ to be $f([c])$. In particular,
$I(a,b) = I([a], [b])$. Homeomorphic manifolds $X$, $Y$ are denoted by $X \cong Y$.
Isotopic submanifolds $c, d$ are denoted by $c \cong d$. If $m \in 
\Cal ML(\Sigma)$,
$I_m$ denotes the geometric intersection number function with
respect to  $m$.
A class in $\Cal CS(\Sigma_{g,r})$ is called 
a \it Fenchel-Nielsen system \rm (resp. an  \it ideal triangulation \rm) 
if it is the isotopy class of $3g+r-3$ (resp. $6g+2r-6$)
many pairwise non-isotopic non-boundary parallel simple loops 
(resp. proper arcs). The numbers $3g+r-3$ and $6g+2r-6$ are maximal.

\it A convention \rm : all surfaces drawn in this paper have the right-hand 
orientation in the front face.

2.1. A multiplicative structure on $\Cal CS(\Sigma)$

Suppose $a$ and $b$ are two arcs in $\Sigma$ intersecting transversely at one
point $P$. Then the resolution of $a \cup b$ at $P$ from $a$ to $b$ is
defined as follows. Take any orientation on $a$ and use the orientation
on $\Sigma$ to determine an orientation on $b$. Then resolve the intersection
according to the orientations. The resolution is independent of the choice 
of the orientation on $a$. See figure 2.

Given two curve systems $a$, $b$ on $\Sigma$ with $|a \cap b| = I(a,b)$, the
multiplication $ab$ is defined to be the disjoint union of simple loops
and arcs obtained by resolving all intersection points from $a$ to $b$. It is 
shown in  \S8 (lemma 8.1) that $ab$ is again a curve system whose isotopy class
depends only on the isotopy classes of $a,$ $b$. Given $\alpha$,$\beta$ $\in$ $\Cal CS(\Sigma)$,
we define $\alpha$$\beta$ =$[ab]$ where $a \in$ $\alpha$, $b \in$ $\beta$ so that $|a \cap b|
= I(a,b)$. The following proposition establishes the basic properties of the
multiplication. See  \S8 for a proof.

Let $\Cal CS_0(\Sigma)$ be the subset of $\Cal CS(\Sigma)$ consisting of
isotopy classes of curve systems which contain no
arcs.

\midspace{0.1cm}
\centerline{\epsfbox{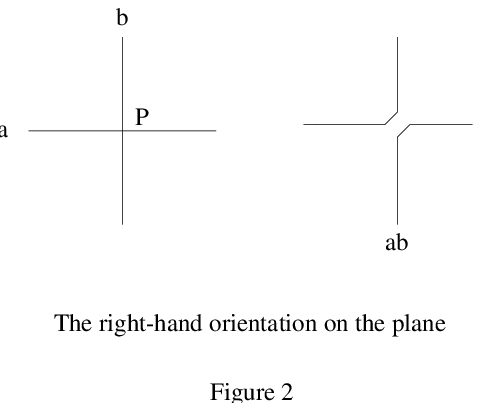}}
\midspace{0.1cm}

{\bf Proposition 2.1.} \it The multiplication $\Cal CS(\Sigma)$$\times$$\Cal CS(\Sigma)$ $\to$ $\Cal CS(\Sigma)$ sends
$\Cal CS_0(\Sigma)$$\times$ $\Cal CS_0(\Sigma)$ to $\Cal CS_0(\Sigma)$ and satisfies the following properties.

(a) It is preserved by the action of the orientation preserving
homeomorphisms.

(b) If $I(\alpha, \beta)$ =0, then $\alpha \beta = \beta \alpha$. 
Conversely, if 
$\alpha \beta = \beta \alpha$ and $\alpha \in$ $\Cal CS_0(\Sigma)$, then $I(\alpha, \beta)$=0. 

(c) If $\alpha$ $\in$ $\Cal CS_0(\Sigma)$, $\beta$ $\in$ $\Cal CS(\Sigma)$, then $I(\alpha, \alpha \beta)=I
(\alpha, \beta \alpha) = I(\alpha, \beta)$ and $\alpha(\beta \alpha)=
(\alpha \beta) \alpha$. If in addition that each component of
$\alpha$ intersects $\beta$, then $\alpha(\beta \alpha) = \beta$.

(d) If $[c_i] \in$ $\Cal CS(\Sigma)$ so that $|c_i \cap c_j| = I(c_i, c_j)$ for $i,j=1,2,3$,
$i \neq j$, $|c_1 \cap c_2 \cap c_3| =0$,
 and there is no contractible region in $\Sigma -(c_1 \cup c_2 \cup c_3)$
bounded by three arcs in $c_1$, $c_2$, $c_3$, then $[c_1]([c_2][c_3])=
([c_1][c_2])[c_3]$.

(e) For any positive integer $k$, $(\alpha ^k \beta^k) =(\alpha \beta)^k$.

(f) If $\alpha$ is the isotopy class of a simple closed curve, then the
positive Dehn twist along $\alpha$ sends $\beta$ to $\alpha^k \beta$ where $k =$ $I(\alpha, \beta)$.
\rm 

It follows from the definition that $I(\alpha, \gamma) + I(\beta, \gamma) 
\geq I(\alpha \beta, \gamma)$.
Furthermore, proposition (c) implies a
stronger result that $I(\alpha \beta, \gamma) + I(\alpha, \gamma) 
\geq I(\beta, \gamma)$ when $\alpha, \beta \in CS_0(\Sigma)$. Indeed,
$I(\alpha \beta, \gamma) + I(\alpha, \gamma) \geq I((\alpha \beta)\alpha, 
\gamma) \geq I(\beta \delta^2, \gamma) \geq I(\beta, \gamma)$
where  $\delta$ consists of components of $\alpha$ which are  disjoint
from $\beta$. 

2.2. The modular relation on $\Cal S(\Sigma_{1,1})$ and $\Cal S(\Sigma_{0,4})$

Call two elements $\alpha$, $\beta$ $\in$  $\Cal S(\Sigma)$ \it orthogonal, \rm denoted by $\alpha \perp \beta$, if
$I(\alpha, \beta)$ = 1; and \it pseudo-orthogonal, \rm denoted by $\alpha \perp_0 \beta$, if $I(\alpha, \beta)$=2 so that
their algebraic intersection number is zero. Suppose 
$\alpha \perp \beta$ or $\alpha \perp_0 \beta$. Take
$a \in$ $\alpha$, $b \in$ $\beta$ so that $|a \cap b|$ = $I(\alpha, \beta)$. Then $N(a \cup b)
\cong \Sigma_{1,1}$ if $\alpha \perp \beta$ and $N(a \cup b) \cong \Sigma_{0,4}$ if $\alpha \perp_0 \beta$.
It follows from the definition that $\alpha \beta \perp \alpha, \beta$ if $\alpha \perp \beta$,
and $\alpha \beta \perp_0 \alpha, \beta$ if $\alpha \perp_0 \beta$.
Thus three distinct elements  $\alpha$, $\beta$, $\gamma$ $\in$ 
$\Cal S'(\Sigma_{1,1})$
(resp. $\Cal S'(\Sigma_{0,4})$) form an ideal triangle if and only if
$\alpha \perp \beta$ (resp. $\alpha \perp_0 \beta$) and $\gamma$ $\in$ $\{\alpha \beta, \beta \alpha\}$. In particular
the distinct ideal triangles in equations (2), (5) in theorem 2 are
$(\alpha_1, \alpha_2, \alpha_1 \alpha_2)$ and $(\alpha_1, \alpha_2,
\alpha_2 \alpha_1)$ where $\alpha_1 \perp \alpha_2$ or $\alpha_1 \perp_0
\alpha_2$ ($(\alpha_1, \alpha_2, \alpha_1 \alpha_2)$ is
positively oriented). 
If $\alpha \perp \beta$ or $\alpha \perp_0 \beta$, we define $\alpha^{-n} \beta = \beta \alpha^n$
for $n \in \bold Z_{>0}$. It follows from proposition 2.1(c) that 
$\alpha^n(\alpha^m \beta) = \alpha^{n+m} \beta$ for $n,m \in \bold Z$.

For $\Sigma = \Sigma_{1,1}$ or $\Sigma_{0,4}$,
we can find an explicit bijection from $\Cal S'(\Sigma)$ to $\hat \bold Q$
as follows. Take $\alpha$, $\beta$ in $\Cal S'(\Sigma)$ so that $\alpha \perp \beta$ or $\alpha \perp_0 \beta$. Then each $\gamma$
in $\Cal S'(\Sigma)$ can be expressed uniquely as $\alpha^p \beta^q$ where
$q \in \bold Z_{\geq 0}$, $p \in \bold Z$ and $p,q$ are relatively prime.
Define $\pi(\gamma) = p/q$ from $\Cal S'(\Sigma)$ to $\hat \bold Q$. 
Then $\pi(\gamma_i) = p_i/q_i$, $i=1,2,$ satisfy $p_1q_2 -p_2 q_1 =
 \pm 1$ if and only if $\gamma_1 \perp \gamma_2$ or $\gamma_1 \perp_0 \gamma_2$.

Given two simple loops $a,b$, we use $a \perp b$ to denote $|a \cap b|$
$ =I(a,b) =1$, and use $a \perp_0 b$ to denote $|a \cap b| = I(a, b) =2$
and $[a] \perp_0 [b]$.

2.3. A gluing lemma

Suppose $\Sigma'$ is an incompressible subsurface of $\Sigma$. We define the 
restriction
map $R (= R_{\Sigma'}^{\Sigma})$ : $\Cal CS(\Sigma)$ $\to$ $\Cal CS(\Sigma')$
 as follows.
Given $\alpha$ in $\Cal CS(\Sigma)$, take $a \in \alpha$ so that $|a \cap \partial \Sigma'|
=I(a, \partial \Sigma')$ and $a \cap \Sigma'$ contains no component parallel
into $\partial \Sigma'$. We define $R(\alpha) =[a |_{\Sigma'}] (:=
 \alpha|_{\Sigma'})$. The restriction map is well defined. Furthermore if
$X \subset Y \subset Z$ are incompressible subsurfaces, then $R^Z_X = R^Y_X
 R^Z_Y$.

{\bf Lemma 2.1} \it (Gluing along a 3-holed sphere) \it Suppose $X$ and
$Y$ are incompressible subsurfaces in $\Sigma$ so that $\Sigma = X \cup Y$ and 
$X \cap Y \cong \Sigma_{0,3}$. Then for any two elements $\alpha_X 
\in \Cal CS(X)$,
$\alpha_Y \in \Cal CS(Y)$ with $\alpha_X|_{X \cap Y} = \alpha_Y|_{X \cap Y}$,
there is a unique element $\alpha \in $ $\Cal CS(\Sigma)$ so that $\alpha|_X = \alpha_X$ and
$\alpha|_Y =\alpha_Y$. \rm

Proof. To show the existence, take $ a_1 \in  \alpha_X $ and $a_2 \in 
\alpha_Y$
so $\alpha_X |_{X \cap Y}$ = $ [a_1|_{X \cap Y}]$, $\alpha_Y |_{X \cap Y}
=[a_2|_{X \cap Y}]$. By the assumption, there is a self-homeomorphism $h_1$
of $X \cap Y$ isotopic to the identity map so that $ h_1(a_1|_{X \cap Y})
= a_2|_{X \cap Y}$. Extend $h_1$ to a self-homeomorphism $h_2$ of $X$ isotopic
to $id_X$. Define a curve system $a$ on $\Sigma$ as follows: $a|_X = h_2(a_1)$, and
$a|_Y = a_2$. Then we have $[a]|_X = \alpha_X$ and $[a]|_Y = \alpha_Y$ by
definition.

To show the uniqueness, suppose $\beta$ $\in$ $\Cal CS(\Sigma)$ so that $\beta|_X = \alpha_X$,
and $\beta|_Y = \alpha_Y$. Take $b \in \beta$ so that $b|_X \in \alpha_X$.
There is a self-homeomorphism $h_3$ of $X$ isotopic to $id_X$ so that
$h_3(b|_X) = a|_X$. By extending $h_3$ to a homeomorphism of $\Sigma$ isotopic
to $id_{\Sigma}$, we may assume that $b|_X = a|_X$. Now since $a|_Y \in
\alpha_Y$ and $b|_X = a|_X$, we obtain $b|_Y \in \alpha_Y$ (due to $\partial Y
\cap int(\Sigma)  \subset int(X)$). Let
$h_4$ be a self-homeomorphism of $Y$ sending $b|_Y$ to $a|_Y$ so that
$h_4 \cong id_Y$ and $h_4|_{\partial Y \cap (\partial(X \cap Y))} = id$.
Extend $h_4$ to a homeomorphism $h_5$ of $\Sigma$ by setting $h_5(x) = x$ for
$x \in X -Y$. Then $h_5 \cong id$ and $h_5(b) = a$. Thus $\alpha$ = $\beta$.
$\square$

\midspace{0.1cm}
\centerline{\epsfbox{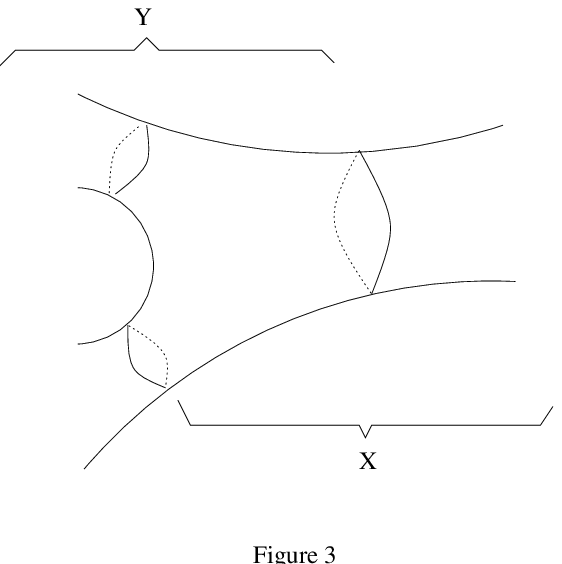}}
\midspace{0.1cm}

\it Remarks 2.1. \rm The lemma also holds for measured laminations. An easy way
to derive it is to use Dehn-Thurston's parametrization of $\Cal ML(\Sigma)$ based
on a Fenchel-Nielsen system of the surface $\Sigma$ so that each component of 
$\partial(X \cap Y)$ is either in the Fenchel-Nielsen system or  is
a boundary component of the surface (see [FLP] or [PH]). 
Given a Fenchel-Nielsen system
$\alpha = \alpha_1 \cup ... \cup \alpha_k$ where $k =3g+r-3$ and two classes
$\beta$, $\gamma$ $\in$ $\Cal CS_0(\Sigma)$ so that $I(\beta, \alpha_i) = 
I(\gamma, \alpha_i)$ for 
all $i$, we can express $\beta = \alpha_i^{n_1} .... \alpha_{k}^{n_k}
\gamma$ where $n_i \in \bold Z$ by the defintion of the multiplication
(recall that $\alpha^{-n} \delta = \delta \alpha^n $ for $n <0$). We call
$(n_1,...,n_k)$ the \it relative Dehn-Thurston twisting
coordinate \rm of $\beta$ with respect to $\gamma$. The twisting coordinates
and the intersection number coordinates $I(\beta, \alpha_i)$ form the Dehn-Thurston
parametrization. Now the proof of the lemma follows easily by comparing the
twisting coordinates at $\partial (X \cap Y)$. 

\it 2.2. \rm For surface with boundary,  Mosher [Mo] has introduced
a parametrization of $\Cal ML(\Sigma)$ using an ideal triangulation
where the coordinates are the intersection numbers.

\S 3. The One-holed Torus 

The goal of this section is to show theorem 2 for $\Sigma_{1,1}$. We restate the result
in terms of the multiplicative structure  as follows.

{\bf Proposition 3.1.} \it A function $f: \Cal S(\Sigma_{1,1}) \to \bold Z_{\geq 0}$
is the geometric intersection number
function $I_{\delta}$ for some $\delta \in \Cal CS(\Sigma_{1,1})$ if and only if
for $\alpha \perp \beta$ and $\gamma = \alpha \beta$,
$$f(\alpha) + f(\beta) + f(\gamma) = \max(2f(\alpha), 2f(\beta), 
2f(\gamma), f(\partial
\Sigma_{1,1})) \tag 1$$ 
$$f(\alpha \beta) + f(\beta \alpha) = \max(2f(\alpha), 2f(\beta), f(\partial
\Sigma_{1,1})) \tag 2$$
$$f(\partial \Sigma_{1,1}) \in 2 \bold Z. \tag 3$$
Furthermore, the characterization of geometric functions  $f:
\Cal S \to \bold R$ is given by equations
(1), (2) above. \rm

\it Remark. \rm The condition $f \geq 0$ in the proposition above is
not necessary. Indeed, equation (1) (also equation (4))
implies $f \geq 0$.  To see this, we note that (1) implies  that $f(\alpha)$,
$f(\beta)$, $f(\gamma)$ satisfy the triangular inequalities that sum of
two is at least the third which in turn shows $f \geq 0$.

Proof. To see the necessity, we double the surface $\Sigma_{1,1}$ to 
obtain $\Sigma_{2,0} =\Sigma_{1,1} \cup_{id_{\partial}} \Sigma_{1,1}$. 
 Then each $\gamma$ $\in$ $\Cal CS(\Sigma_{1,1})$ 
corresponds to $\hat \gamma \in$ $\Cal CS(\Sigma_{2,0})$ whose restriction
to both summands $\Sigma_{1,1}$ are $\gamma$. The curve system  
$\hat \gamma$ has no boundary. Let $d_i$
be a sequence of a hyperbolic metrics on $\Sigma_{2,0}$ which pinch to 
$\hat \gamma$,
i.e., there is a sequence $\lambda_i \in \bold R_{>0}$ so that 
$\lim_{i} \lambda_i l_{d_i}(\alpha) = I_{\hat \gamma}(\alpha)$
for all $\alpha \in$ $\Cal S(\Sigma_{2,0})$ where $l_{d_i}(\alpha)$ is the length of the
$d_i$-geodesic in the class $\alpha$. Let $t_i = 2 \cosh l_{d_i}/2$. It is
shown in  [FK], [Ke] and [Lu1] that for $\alpha \perp \beta$ 
in $\Cal S(\Sigma_{1 ,1})$ ($\subset \Cal S(\Sigma_{2,0}))$,
one has the following identities: 
$t_i(\alpha) t_i(\beta) t_i(\alpha \beta) = t_i^2(\alpha) + t_i^2(\beta)
+ t_i^2(\alpha \beta) + t_i(\partial \Sigma_{1,1}) -2$ and $t_i(\alpha \beta) t_i(\beta \alpha)
=t_i^2(\alpha) + t_i^2( \beta) + t_i(\partial \Sigma_{1,1}) -2$. Now, for
$\alpha$ $\in$ $\Cal S(\Sigma_{1,1})$, we have $I_{\hat \gamma}(\alpha) = I_{\gamma}(\alpha)$. Let
$i$ tend to infinity. The equations for $t_i$ degenerate to the equations (1),
(2) in the proposition. The equation (3) is evident.

\it Remark 3.1. \rm To derive equation (1) directly from the trace
identity $tr(A) tr(B)$ \newline
$ tr(AB) = tr^2(A) + tr^2(B) + tr^2(AB)
- tr([A,B]) -2$ where $A, B \in SL(2, \bold R)$,  we assume that
$A, B, AB$ correspond to three simple closed geodesics forming an ideal 
triangle in $\Cal S$. Then
 $tr(A)tr(B)tr(AB) > 0$ and $tr([A,B]) < 0$ (see [GiM] for instance).
In  particular, we obtain
$|tr(A) ||tr(B)||tr(AB)| = tr^2(A) + tr^2 (B) + tr^2(AB) + |tr([A,B])| -2$.
The degeneration of it becomes $f(A)
+ f(B) + f(AB) =\max(2f(A), 2f(B), 2f(AB),f([A,B]))$ which is equation (1).

\smallskip
To show that the  conditions are also sufficient, we begin with a  function
$f : \Cal S \to \bold Z_{\geq 0}$ satisfying equations (1),(2),(3). 
By the structure  of the modular relation, we conclude that $f$ 
is determined
by its values on $\{\alpha, \beta, \alpha \beta, \partial \Sigma_{1,1}\}$
for $\alpha \perp \beta$. Thus it suffices to construct $\delta \in $ $\Cal 
CS(\Sigma_{1,1})$ so that $f$ and
$I_{\delta}$ have the same values at the four-element set above. 

We consider two cases: $\min\{f(\alpha): \alpha \in \Cal S'(\Sigma_{1,1})\} =0
$, or $>0$.

Case 1. There is  $\alpha$ $\in$ $\Cal S'(\Sigma)$ so that $f(\alpha) =0$. If $\beta \perp
\alpha$ and $\gamma = \alpha \beta$, then $f(\beta) = f(\gamma)$. Indeed,
by equation (1), $f(\beta) + f(\gamma) =  \max(2f(\beta), 2f(\gamma), 
f(\partial \Sigma_{1,1}))$ $\geq \max(2f(\beta), 2f(\gamma))$. Thus $f(\beta) = f(\gamma)$.
In particular, $f(\beta) \geq \frac{1}{2} f(\partial \Sigma_{1,1})$. We construct
the curve system $\delta$ as follows. Let $\Sigma'$  $=\Sigma_{1,1} - int(N(a))$
where $a \in \alpha$. Then $\Sigma'$ $\cong \Sigma_{0,3}$. Curve systems on $\Sigma_{0,3}$
with $\partial \Sigma_{0,3} = b_1 \cup b_2 \cup b_3$ are well 
understood. Namely, 
$\Cal S(\Sigma_{0,3})$ $=\{b_1, b_2, b_3\}$ and each $\delta \in$ $\Cal 
CS(\Sigma_{0,3})$ is uniquely determined by
$\pi(\delta) =(I_{b_1}(\delta), I_{b_2}(\delta), I_{b_3}(\delta))$. Furthermore,
each triple of non-negative integers whose sum is even is of the form
$\pi(\delta)$ and $\pi(\delta \delta') = \pi(\delta) + \pi(\delta')$.
Let $\delta' \in \Cal CS(\Sigma')$ $( \subset \Cal SC(\Sigma))$
so that $I(\delta', \partial \Sigma_{1,1})
=f(\partial \Sigma_{1,1})$ and $I(\delta', \alpha) =0$. Let $\delta =
\delta' \alpha^k$ in $\Cal CS(\Sigma)$ where $k= f(\beta) -\frac{1}{2} f(\partial \Sigma_{1,1})$.
Then $I_{\delta}$ and $f$ have the same values at $\{\alpha, \beta, \gamma,
\partial \Sigma_{1,1}\}$ by the construction. Thus $f =I_{\delta}$.

Case 2. Suppose $\min\{f(\alpha): \alpha \in \Cal S'(\Sigma_{1,1})$\} $>0$.
Let $\alpha \perp \beta$ be the classes so that $f(\alpha) + f(\beta) +f(\alpha \beta)
 = \min\{ f(\alpha') + f(\beta') + f(\alpha' \beta'): \alpha' \perp \beta'\}.$
We claim that $f(\alpha) + f(\beta) + f(\alpha \beta) = 
f(\partial \Sigma_{1,1})$.
To see this, let $\gamma =\alpha \beta$ and we assume without loss of generality
that $f(\alpha) \geq f(\beta) \geq f(\gamma) > 0$ (since $\{\alpha, \beta,
\gamma\}$ is symmetric). Suppose 
the claim is false. Then equation (1) shows
that $f(\alpha) + f(\beta)
+f(\gamma) > f(\partial \Sigma_{1,1})$. Furthermore, equation (1) shows that
$f(\alpha) = f(\beta) + f(\gamma)$ $> \max(f(\beta), f(\gamma))$.  It follows
$f(\partial \Sigma_{1,1}) < f(\alpha) + f(\beta) + f(\gamma) = 2f(\alpha)$. 
Consider equation (2) for $\alpha$ (=$\beta \gamma$) and $\alpha' 
(= \gamma \beta$). We obtain
$f(\alpha) + f(\alpha') = \max(2f(\beta), 2f(\gamma), f(\partial \Sigma_{1,1}))$
$< 2f(\alpha)$. Thus $f(\alpha') < f(\alpha)$ which contradicts the choice of
$\{\alpha, \beta, \gamma\}$.

Now equation (1) shows  that $f(\alpha), f(\beta), f(\gamma)$ satisfy the triangular
inequalities (sum of two is no less than the third) and their  sum is an 
even  number. Thus there exist  integers $x,y,z \in \bold Z_{\geq 0}$ so that
$f(\alpha) = y+z$, $f(\beta) = z + x$, and $f(\gamma) = x + y$. Let
$\alpha_1 \beta_1 \gamma_1$ in $\Cal CS(\Sigma)$ be the ideal triangulation so that
$I(\alpha, \alpha_1) = I(\beta, \beta_1) = I(\gamma, \gamma_1) =0$ (see 
figure 4). Define $\delta = \alpha_1^x \beta_1^y \gamma_1^z$. Then
$f = I_{\delta}$ on the four element set $\{\alpha, \beta, \gamma, 
\partial \Sigma_{1,1}\}$ . Thus $f =I_{\delta}$.

\midspace{0.1cm}
\centerline{\epsfbox{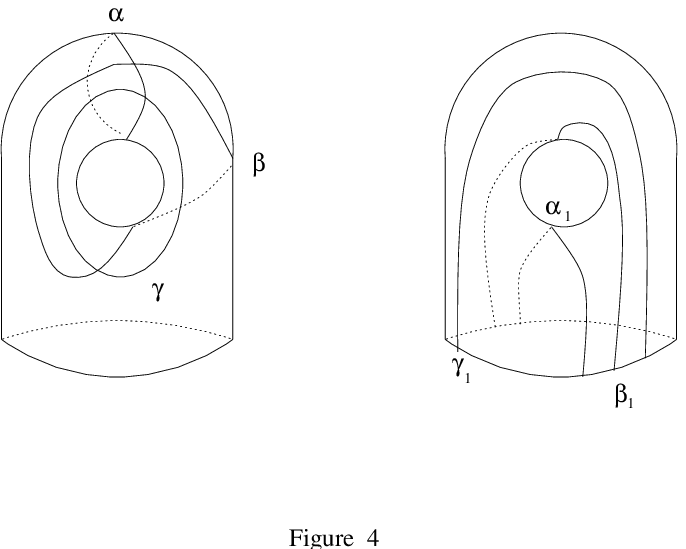}}
\midspace{0.1cm}

\bigskip
\bigskip

To show part (c), we need the following lemmas.

{\bf Lemma 3.1.} \it The equation $x + a = \max(2x, x+b, c)$ has solutions
in $x$ over $\bold R$ if and only if $a \geq \max(b, c/2)$. If it has 
solutions, then the set of all solutions is given by (i) $\{c-a, a\}$  in 
the case $a > b$, and  by (ii) the closed interval $[c-a, a]$ in the
case $a=b$. In particular, we have
(a) if $x_1$ is a solution, then $c-x_1$ is also a  solution;
(b) if $x_1$ and $x_2$ are solutions so that $x_1 + x_2 = c$, then
$\max(x_1, x_2, b) = a$. \rm

Proof. If $a \geq \max(b,c/2)$, then $x=a$ is a solution. If $x'$ is a solution,
then since $x'+a \geq  x'+b$, we have $a \geq b$. Also $x'+a \geq 2x'
$ and $x'+a \geq c$. Thus $ a \geq x' \geq c-a$. This shows $a \geq c/2$, i.e.,
$a \geq \max( b, c/2)$. If $a >b$, then the equation becomes $x+a =\max(2x,
c)$ with $a \geq c/2$. Thus the solutions are $\{c-a, a\}$. If $a=b$, then
one checks easily that all solutions are points in $[c-a, a]$.
$\square$

{\bf Lemma 3.2.} \it Suppose $x_1, x_2,$$ x_3,$$ x_4$$ \in$$ \bold Z_{\geq 0}$ 
so that $x_1 $$+ x_2$$ + x_3$$=\max$($2x_1,$$2x_2,$$2x_3$,\newline
$x_4$). 
Then there is a function
$g : \Cal S(\Sigma_{1,1}) \to \bold Z$ satisfying equations (1),(2) and an
ideal triangle $(\alpha_1, \alpha_2, \alpha_3)$ in $\Cal S'(\Sigma_{1,1})$
so that $g(\alpha_i) = x_i$, $i=1,2,3$, and $g(\partial \Sigma_{1,1}) =x_4$.
\rm

Proof. Take any ideal triangle  $(\alpha_1, \alpha_2, \alpha_3)$. We define $g$
on $\alpha_i$  and $\partial \Sigma_{1,1}$ as required. We now extend
$g$ through the neighboring ideal triangles by using equation (2).  
Thus, we need to verify  that the equation (1)  for $g$ on the 
neighboring ideal triangles  still holds. Take a neighboring
ideal triangle, say $(\alpha_1, \alpha_2, \alpha_3')$.  
Define $g(\alpha_3') =
x_3'$ where $x_3' = \max(2x_1, 2x_2, x_4) - x_3$. We first note that $x'_3
\geq 0$ since $x_i+ x_j \geq x_k$ for $\{i,j,k\} = \{1,2,3\}$ by the
given condition on $x_i's$. Next, consider
$x_1 + x_2 + x_3  = \max(2x_1, 2x_2, 2x_3, x_4)$ as an equation in $x_3$.
Then it is of the form $x+ x_1 + x_2 = \max(2x, x_3 + x_3')$. 
By lemma 3.1(a), $x_3'$ satisfies the equation in $x$, i.e.,
equation (1) holds for $g$ on the neighboring ideal triangles.
$\square$

We now show that equations (1), (2) characterize the geometric functions. Evidently,
any geometric functions satisfies the equations (1), (2). Conversely,
suppose that $f$ is a solution to equations (1), (2). Fix an ideal triangle
 $(\alpha_1, \alpha_2, \alpha_3)$ in $\Cal S'$. Note that the 
rational solutions
of the equation $x_1 + x_2 + x_3  = \max(2x_1, 2x_2, 2x_3, x_4)$  are
dense in the solutions over $\bold R_{\geq 0}$.  By lemma 3.2, there is
a sequence of functions $g_n$ from $\Cal S$ to $2 \bold Z_{\geq 0}$ solving
equations (1), (2) and a sequence of numbers $k_n \in \bold Q$ so that
$\lim_{n} k_n g_n(x) = f(x)$ for $x \in \{\alpha_1, \alpha_2, \alpha_3,
\partial \Sigma\}$. By equation (2), we have $\lim_n k_n g_n(x) = f(x)$ for all
$x \in$ $\Cal S(\Sigma)$. On the other hand, we have $g_n = I_{\delta_n}$ for some
$\delta_n \in \Cal S(\Sigma)$ by the result for curve systems. Thus
$f =I_{m}$ where $m = \lim_n k_n \delta_n$ $\in$ $\Cal ML(\Sigma)$ by definition. $\square$

\S 4. The Four-holed Sphere

The goal  of this section is to show theorem 2 for the surface
$\Sigma_{0,4}$. The basic ideal of the proof is the same as  in \S3. But
the proof is considerably longer and more complicated due to the
existence of eight non-homeomorphic ideal triangulations of the
four-holed sphere.
We restate the theorem in terms of the multiplicative structure below.

{\bf Proposition 4.1.} \it For surface $\Sigma_{0,4}$ with $\partial \Sigma_{0,4}
= b_1 \cup b_2 \cup b_3 \cup b_4$, a function $f :\Cal S(\Sigma_{0,4})$ 
$\to \bold Z_{\geq 0}$  is the geometric intersection number function $I_{\delta}$
for some $\delta \in$ $\Cal CS(\Sigma)$ if and only if for $\alpha_1 \perp_0 \alpha_2$
with $\alpha_3 = \alpha_1 \alpha_2$ so that $(\alpha_i, b_s, b_r)$
bounds a $\Sigma_{0,3}$ in $\Sigma_{0,4}$,
$$\sum_{i=1}^3 f(\alpha_i) = \max_{1 \leq i \leq 3;1 \leq s \leq 4}
(2f(\alpha_i), 2f(b_s), \sum_{s=1}^4 f(b_s), f(\alpha_i)+ f(b_s) + f(b_r))
\tag 4$$
$$f(\alpha_1 \alpha_2) + f(\alpha_2 \alpha_1) =\max_{1 \leq i \leq 2;1 \leq s \leq 4}
(2f(\alpha_i), 2f(b_s), \sum_{s=1}^4 f(b_s), f(\alpha_i)+ f(b_s) + f(b_r))
\tag 5$$
$$f(\alpha_i) + f(b_s) + f(b_r) \in 2 \bold Z \tag 6$$
Furthermore, geometric functions on $\Cal S(\Sigma_{0,4})$ are characterized by the equations
(4),(5). \rm

Proof. The necessity of the equations (4),(5) follows from the same argument
as in \S3 using the degenerations of the trace relations for geodesic
length functions. To be more precise, it is shown in [Lu1] that for
any hyperbolic metric $d$ on $\Sigma_{0,4}$ with geodesic boundary or cusp ends,
then $t(\alpha) = 2 \cosh l_d(\alpha)/2$ satisfies:
$$t(\alpha_1) t(\alpha_2) t(\alpha_3) + 4 =  \sum_{i=1}^3 t^2(\alpha_i) + 
\sum_{s=1}^4t^2(b_s ) + \prod_{s=1}^4 t(b_s) +   \frac{1}{2}
\sum_{i=1}^3 \sum_{s=1}^4 t(\alpha_i) t(b_s) t(b_r)$$ 
and
$$t(\alpha_1 \alpha_2)t(\alpha_2 \alpha_1) = \sum_{i=1}^2 t^2(\alpha_i)
 + \sum_{s=1}^4t^2(b_s ) + {\dsize \prod_{s=1}^4} t(b_s)
+ {\dsize  \frac{1}{2} \sum_{i=1}^2 \sum_{s=1}^4} t(\alpha_i) t(b_s) t(b_r) -4$$ 
Now the degenerations of the above two equations are equations (4), (5). The
equation (6) holds for curve systems clearly.

To show that the  conditions are also sufficient, we begin with a  function
$f : \Cal S \to \bold Z_{\geq 0}$ satisfying equations (4),(5),(6). By the structure  of the
modular relation, we conclude that $f$ is determined by its restriction
on $\{\alpha, \beta, \alpha \beta, b_1,..., b_4\}$ for $\alpha \perp_0 \beta$.
Thus it suffices to construct $\delta$ $\in$ $\Cal CS(\Sigma)$ so that $f$ and $I_{\delta}$ have the
same values on the seven-element set  $\{\alpha, \beta, \alpha \beta, b_1,..., b_4\}$.

Note that equation (6) implies both $\sum_{i=1}^4 f( b_i)$ and $\sum_{i=1}^3 f(\alpha_i)$ are even numbers.

We shall consider two cases: $\min\{f(\alpha): 
\alpha \in \Cal S'(\Sigma_{0,4})\} =0$ or $>0$.

Case 1. Suppose $f(\alpha) = 0$ for some  $\alpha \in \Cal S'(\Sigma_{0,4})$. 
Choose $\beta$ so that $\beta \perp_0 \alpha$ and $\gamma =
\alpha \beta$. Then $f(\beta) = f(\gamma)$ due to equation (4) that
$f(\beta) + f(\gamma) = \max(2f(\beta), 2f(\gamma), *) \geq \max(2f(\beta),
2f(\gamma))$. Assume without loss of generality that $(\alpha, b_1, b_2)$,
$(\beta, b_1, b_3)$ bound $\Sigma_{0,3}$ in $\Sigma_{0,4}$. Construct a curve
system $\delta' \in $ $\Cal CS(\Sigma_{0,4})$ so that $I(\delta', \alpha) =0$, $I(\delta', b_i)
= f(b_i)$. The existence of $\delta'$ is due to the classification of curve
systems on $\Sigma_{0,3}$ and the equation (6) that $f(b_1) + f(b_2), f(b_3) + f(b_4)$ are
even numbers. Let $k = \frac{1}{2}( f(\beta) - \max(f(b_1), f(b_2)) -
\max(f(b_3), f(b_4)))$. Then equation (4) for $\alpha \perp_0 \beta$ shows that
$k \geq 0$ and equation (6) shows that $k \in \bold Z$. Let $\delta$ $=\delta' \alpha^k
\in \Cal CS(\Sigma_{0,4})$. Then $I_{\delta}$ and $f$ have the same values on the
set $\{\alpha, \beta, \alpha \beta, b_1,..., b_4\}$.

\midspace{0.1cm}
\centerline{\epsfbox{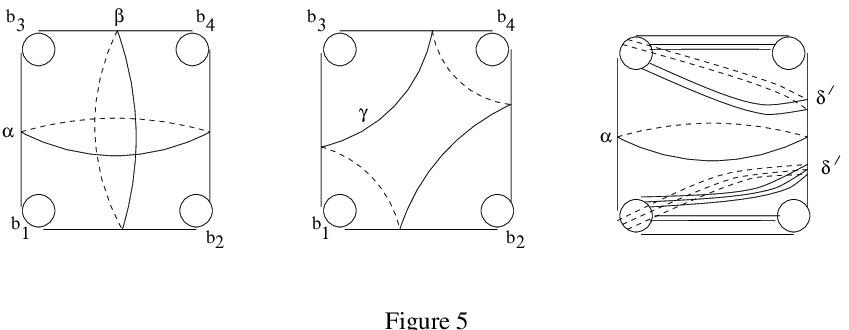}}
\midspace{0.1cm}

Case 2. Assume $f(\alpha) \geq 1$ for all $\alpha \in \Cal S'(\Sigma_{0,4})$.
Let $\{\alpha, \beta, \gamma\}$
be an ideal triangle in $\Cal S(\Sigma)$ so that $f(\alpha) + f(\beta) + f(\gamma)$
achieves the minimal values among all such triples. Assume without loss of generality that $(\alpha, b_1, b_2)$,
$(\beta, b_1, b_3)$ bound $\Sigma_{0,3}$ in $\Sigma_{0,4}$ and that $f(\alpha) \geq
f(\beta) \geq f(\gamma)$. We claim that $f(\alpha) + f(\beta) + f(\gamma)$
$=A$ where $A = \max_{1 \leq s \leq 4}(2f(b_s), \Sigma_{s=1}^4 f(b_s),
f(\alpha)+ f(b_1) +f(b_2), f(\alpha) + f(b_3) + f(b_4), f(\beta) + f(b_1) + f(b_3), f(\beta) + f(b_2) + f(b_4)$, $f(\gamma) + f(b_1) + f(b_4)$,
$f(\gamma) + f(b_2) + f(b_3))$. Indeed, if otherwise,
by equation (4) that  $f(\alpha) + f(\beta) + f(\gamma)$ $= \max(2f(\alpha), 2f(\beta), 2f(\gamma), A)$, we obtain  $f(\alpha) + f(\beta) + f(\gamma)$ $> A$
and  $f(\alpha) = f(\beta) + f(\gamma)$. In particular, $f(\alpha) > f(\beta),
f(\gamma)$, and $2f(\alpha) >  A$. 
Applying equation (5) to $\alpha$, $\alpha'$ where
$\{ \alpha, \alpha'\}= \{\beta \gamma, \gamma \beta$\}, 
we obtain $f(\alpha) + f(\alpha') = \max( 2f(\beta), 2f(\gamma) ,A')$
where $A' \leq A < 2f(\alpha)$. Thus $f(\alpha) + f(\alpha') < 2f(\alpha)$,
i.e., $f(\alpha') < f(\alpha)$. This contradicts the choice of
 $(\alpha, \beta, \gamma)$.

We now construct $\delta$ $\in$ $\Cal CS(\Sigma_{0,4})$ so that $f$ and $I_{\delta}$ have the same
values on \{$\alpha, \beta, \gamma, b_1,..., b_4\}$ under the  assumption
that  $f(\alpha) + f(\beta) + f(\gamma) = A$. For simplicity, we still
assume that $(\alpha, b_1, b_2)$ and $(\beta, b_1, b_3)$ bound
$\Sigma_{0,3}$ but do not assume that $f(\alpha) \geq f(\beta) \geq f(\gamma)$.

By symmetry, since  $f(\alpha) + f(\beta) + f(\gamma)= A$, it suffices to
consider the following three subcases: (2.1)  $f(\alpha) + f(\beta) + f(\gamma)$$= \Sigma_{s=1}^4 f(b_s)$; (2.2)  $f(\alpha) + f(\beta) + f(\gamma)$
$= 2f(b_1)$; and (2.3)  $f(\alpha) + f(\beta) + f(\gamma) = f(\alpha)
+ f(b_1) + f(b_2)$. The corresponding  curve system $\delta$ in $\Cal CS(\Sigma_{0,4})$ will be constructed as follows. First, we construct an ideal triangulation
$\tau = \tau_1 ... \tau_6$ of $\Sigma_{0,4}$. Then the curve system 
$\delta$ is taken
to be of the form $\tau_1^{x_1}... \tau_6^{x_6}$, $x_i \in \bold Z_{\geq 0}$.

Case (2.1).  $f(\alpha) + f(\beta) + f(\gamma) = \sum_{s=1}^4 f(b_s)$. 
The ideal triangulation $\tau$ is as shown in figure 6 where the locations of
$\alpha$, $\beta$, $\gamma$ are indicated. The conditions that $f$ and
$I_{\delta}$ have the same values on $\{ \alpha, \beta, \gamma, b_1,..., b_4\}$
are given by the following systems of linear equations in $x_i$.

\midspace{0.1cm}
\centerline{\epsfbox{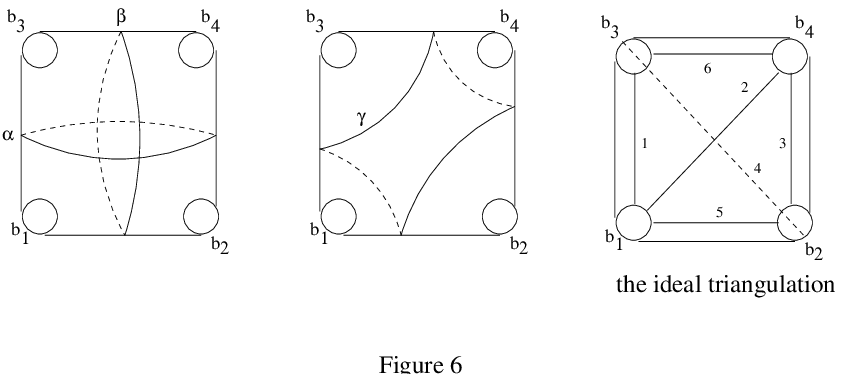}}
\midspace{0.1cm}

$$
\gather
x_1 + x_2 + x_5  = f(b_1)\\
x_3 + x_4 + x_5 = f(b_2)\\
x_1 + x_4 + x_6  = f(b_3)\\
x_2 + x_3 + x_6  = f( b_4)\\
x_1 + x_2 + x_3 + x_4  = f(\alpha)\\
x_2 + x_4 + x_5 + x_6  = f(\beta)\\
x_1 + x_3 + x_5 + x_6  = f(\gamma)\\
\endgather
$$
Note that  $f(\alpha) + f(\beta) + f(\gamma) = \Sigma_{s=1}^4 f(b_s)$ 
is a consequence of the equations above. Thus, it is essentially a systems
of six equations in six variables. The solution is
$$\gather
x_1 = (f(b_1) + f(b_3) - f(\beta))/2\\
x_2 = (f(b_1) + f(b_4) - f(\gamma))/2\\
x_3 = (f(b_2) + f(b_4) - f(\beta))/2\\
x_4 = (f(b_2) + f(b_3) -f(\gamma))/2\\
x_5 = (f(b_1) + f(b_2) - f(\alpha))/2\\
x_6 = (f(b_3) + f(b_4) - f(\alpha))/2\endgather
$$
It remains to show that $x_i \in \bold Z_{\geq 0}$. First of all $x_i \in  \bold
Z$ due to equation (6). To see $x_i \geq 0$, say $x_1 \geq 0$, for definiteness,
we use equation (4) that  $f(\alpha) + f(\beta) + f(\gamma) \geq f(\beta) + f(b_2)
+ f(b_4)$. But  $f(\alpha) + f(\beta) + f(\gamma) = \Sigma_{s=1}^4 f(b_s)$.
Thus, $f(b_1) + f(b_3) \geq f(\beta)$, i.e., $x_1 \geq 0$. The proof of the
rest of the cases $x_i \geq 0$ is similar.  (The solutions $x_i$  are found
as follows: $x_1$ is the number of arcs joining $b_1$, $b_3$ in the
3-holed sphere $\Sigma_{0,3}$ bounded by $b_1$, $b_3$, $\beta$, etc.).

Case (2.2).  $f(\alpha) + f(\beta) + f(\gamma) = 2f(b_1)$. The curve
system $\delta$ is based on the ideal triangulation $\tau$ as shown in figure 7.
We obtain the following system of linear equations in $x_i$

\midspace{0.1cm}
\centerline{\epsfbox{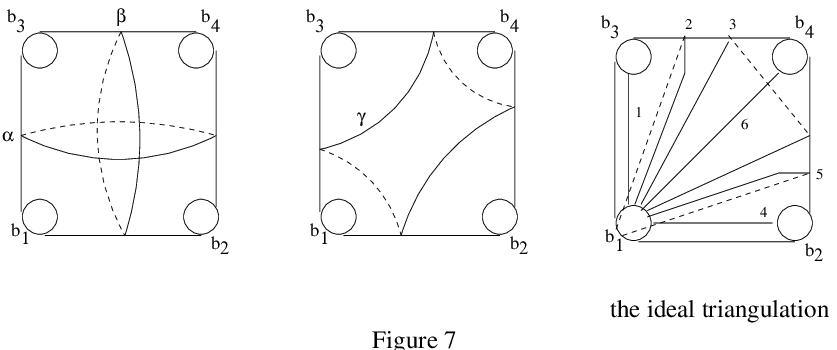}}
\midspace{0.1cm}

$$\gather
x_1 + 2x_2 + 2 x_3 + x_4 + 2x_5 + x_6 = f(b_1)\\
x_4 = f(b_2)\\
x_1 = f(b_3)\\
x_6 = f(b_4)\\
x_1 + 2x_2 + 2x_3 + x_6 = f(\alpha)\\
2x_3 + x_4 + 2x_5 + x_6 = f(\beta)\\
x_1 + 2x_2 + x_4 + 2x_5 = f(\gamma)\endgather
$$
The solution is,
$$ \gather
x_1 = f(b_3)\\
x_2 = (f(b_1) - f(b_3) -f(\beta))/2\\
x_3 = (f(b_1) -f(b_4) -f(\gamma))/2\\
x_4 = f(b_2)\\
x_5 = (f(b_1) - f(b_2) -f(\alpha))/2\\
x_6 = f(b_4)\endgather
$$
To see that $x_i \in  \bold Z_{\geq 0}$, we note that $x_i \in
 \bold Z$ by equation (6). To
show $x_i \geq 0$, say $x_2 \geq 0$, we use equation (4) and the assumption that
 $f(\alpha) + f(\beta) + f(\gamma) = 2f(b_1)$. Thus $2f(b_1) \geq 
f(\beta) + f(b_3) + f(b_1)$, i.e., $x_2 \geq 0$. By symmetry, $x_3, x_5
\geq 0$. 

Case (2.3).  $f(\alpha) + f(\beta) + f(\gamma) =f(\alpha) + f(b_1) + f(b_2)$,
i.e., $f(\beta) + f(\gamma) = f(b_1) + f(b_2)$. We first observe 
that many inequalities
follow from the assumption. To simplify the notions, we use $\Delta
=\{(a_1, a_2, a_3) \in \bold R_{\geq 0}: a_i + a_j \geq a_k,
i \neq j \neq k \neq i\}$. For instance, equations (1),(4) show that
 $(f(\alpha_1), f(\alpha_2), f(\alpha_3)) \in \Delta$. 

{\bf Lemma 4.1.} \it Under the assumption  $f(\alpha) + f(\beta) + f(\gamma)
 =f(\alpha) + f(b_1) + f(b_2)$,
we have

(a) $(f(\alpha), f(b_1), f(b_2)) \in \Delta$;

(b) $f(\alpha) \geq f(b_3) + f(b_4)$;

(c) $f(\beta) + f(b_1) \geq f(b_3)$, and $f(\beta) + f(b_2) \geq f(b_4)$;

(d) $f(\gamma) + f(b_1) \geq f(b_4)$ and $f(\gamma) + f(b_2) \geq f(b_3)$.
\rm

Proof. To see (a), since  $(f(\alpha),  f(\beta),  f(\gamma)) \in \Delta$,
thus $f(\alpha) \leq f(\beta) + f(\gamma) = f(b_1) + f(b_2)$. On
the other hand, equation (4) shows that  $f(\alpha) + f(\beta) + f(\gamma) 
\geq 2f(b_i)$, for $i=1,2$. Thus $f(\alpha) + f(b_i) \geq f(b_j)$ for
$\{i,j\} = \{1,2\}$. To see (b), we use  
 $f(\alpha) + f(\beta) + f(\gamma) \geq \Sigma_{s=1}^4 f(b_s)$ and the
assumption. To see   $f(\beta) + f(b_1) \geq f(b_3)$ in part (c), 
we use  $f(\alpha) + f(\beta) + f(\gamma)
\geq f(\gamma) + f(b_2) + f(b_3)$ (by equation (4)). Now  $f(\alpha) + f(\beta) + f(\gamma) \leq f(\beta) + f(\gamma) + f(\beta) + f(\gamma) = f(\beta) + f(\gamma)
+ f(b_1) +f(b_2)$. Thus the result follows. The rest of the inequalities
in (c),(d) are proved by the same argument. $\square$

To construct the curve system $\delta$, we shall consider nine subcases due to
the different situations: $(f(\beta), f(b_i), f(b_j)) \in \Delta$,
$f(\beta) + f(b_i) \geq f(b_j)$ for $(i,j) \in \{(1,3), (3,1), (2,4), (4,2)\}$.
The nine subcases are listed in figure 8. The (i,j)-th subcase 
corresponds to the i-th row and j-th column in 
figure 8. Due to symmetry,  the (i,j)-th subcase and the (j,i)-th subcase are
essentially the same. We shall consider six subcases: (1,1), (1,2), (1,3),
(2,2), (2,3), (3,3). The corresponding ideal triangulations and
the system of linear equations are listed below.

\midspace{0.1cm}
\centerline{\epsfbox{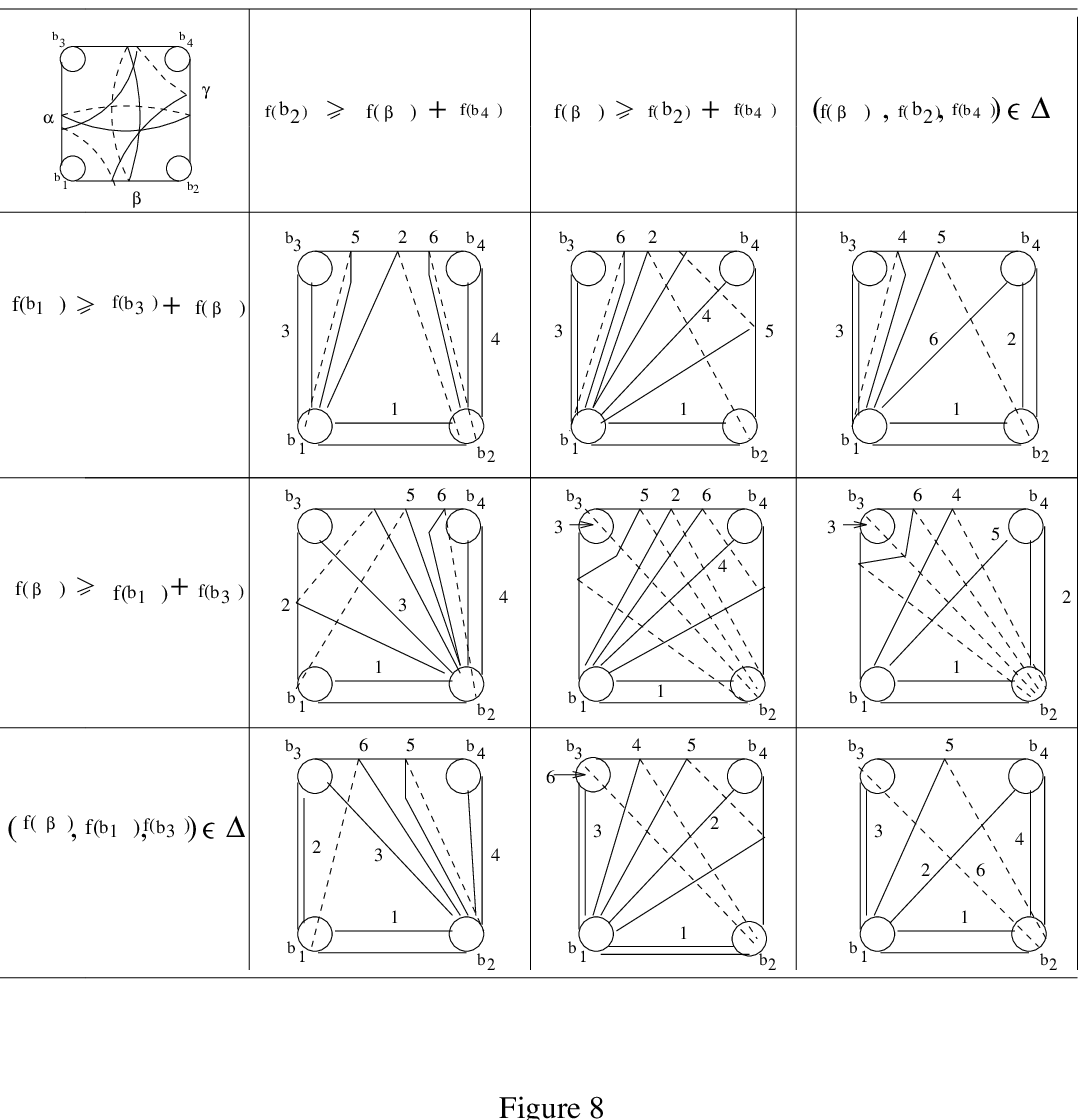}}
\midspace{0.1cm}

Subcase (1,1). $f(b_1) \geq f(\beta) + f(b_3)$ and $f(b_2) \geq
f(\beta) + f(b_4)$.
$$\gather
x_1 + x_ 2+ x_3 + 2x_5 = f(b_1)\\
x_1 + x_2 + x_4 + 2 x_6 = f(b_2)\\
x_3 = f(b_3)\\
x_4 = f(b_4)\\
2x_2 + x_3 + x_4 + 2x_5 + 2x_6 = f(\alpha)\\
x_1 + x_2 = f(\beta)\\
x_1 + x_2 + x_3 + x_4 + 2x_5 + 2x_6 = f(\gamma)\endgather
$$
The solution is,
$$\gather
x_1 = (f(b_1) + f(b_2) -f(\alpha))/2\\
x_2 = (f(\alpha) + f(\beta) -f(\gamma))/2\\
x_3 = f(b_3)\\
x_4 = f(b_4)\\
x_5 = (f(b_1) -f(\beta) -f(b_3))/2\\
x_6 = (f(b_2) -f(\beta) -f(b_4))/2\endgather
$$
The solutions $x_i$ are in $\bold Z_{\geq 0}$ by lemma 4.1, equation (6) and
the assumption ($x_5, x_6 \geq 0$).

Subcase (1.2). $f(b_1) \geq f(\beta) + f(b_3)$ and $f(\beta) \geq f(b_2)+ f(b_4)$.
$$\gather
x_1 + x_2 + x_3 + x_4 + 2x_5 + 2x_6 = f(b_1)\\
x_1 + x_2 = f(b_2)\\
x_3 = f(b_3)\\
x_4 = f(b_4)\\
2x_2 + x_3 + x_4 + 2x_5 + 2x_6 = f(\alpha)\\
x_1 + x_2  + x_4 + 2x_5 = f(\beta)\\
x_1 + x_2 + x_3  + 2x_6 = f(\gamma)\endgather
$$
The solution is,
$$\gather
x_1 = (f(b_1) + f(b_2) -f(\alpha))/2\\
x_2= (f(b_2) + f(\alpha) -f(b_1))/2\\
x_3 = f(b_3)\\
x_4 = f(b_4)\\
x_5 = (f(\beta) -f(b_2) -f(b_4))/2\\
x_6 = (f(b_1) -f(b_3) -f(\beta))/2\endgather
$$
The solutions are in $\bold Z_{\geq 0}$ by lemma 4, equation (6) and the
assumption.

Subcase (1.3). $f(b_1) \geq f(\beta) + f(b_3)$ and $(f(\beta), f(b_2), f(b_4))
\in \Delta$.
$$\gather
x_1 + x_3 + 2x_4 + x_5 + x_6 = f(b_1)\\
x_1 + x_2 + x_5 = f(b_2)\\
x_3 = f(b_3)\\
x_2 + x_6 = f(b_4)\\
x_2 + x_3 + 2x_4 + 2x_5 + x_6 = f(\alpha)\\
x_1 + x_5 + x_6 = f(\beta)\\
x_1 + x_2 + x_3 + 2x_4 + x_5 = f(\gamma)\endgather
$$
The solution is,
$$\gather
x_1 = (f(b_1) + f(b_2) -f(\alpha))/2\\
x_2 =( f(b_2) + f(b_4) -f(\beta))/2\\
x_3 = f(b_3)\\
x_4 = (f(b_1) -f(b_3) -f(\beta))/2\\
x_5 = (f(\alpha) + f(\beta) - f(b_1) - f(b_4))/2\\
x_6 =( f(b_4) + f(\beta) - f(b_2))/2\endgather
$$
By the same argument as in the previous cases, all $x_i$ except 
possibly $x_5$ are in $\bold Z_{\geq 0}$.  It remains to show that $x_5 \in \bold Z_{\geq 0}$. Indeed,
$f(\alpha) + f(\beta) -f(b_1) -f(b_4) = ( f(\alpha) + f(\beta) + f(\gamma))
-(f(\gamma) + f(b_1) + f(b_4))$. Thus, by equations (4), (6), $x_5
\in \bold Z_{\geq 0}$.

Subcase (2.2). $f(\beta) \geq f(b_1) + f(b_3)$ and $f(\beta) \geq f(b_2)+ 
f(b_4))$.
$$\gather
x_1 + x_2 + x_4 + 2x_6 = f(b_1)\\
x_1 + x_2 + x_3 + 2x_5 = f(b_2)\\
x_3 = f(b_3)\\
x_4 = f(b_4)\\
2x_2 + x_3 +  x_4 + 2 x_5 + 2x_6 = f(\alpha)\\
x_1 + x_2 + x_3 + x_4 + 2x_5 + 2x_6 = f(\beta)\\
x_1 + x_2  = f(\gamma)\endgather
$$
The solution is
$$\gather
x_1 = (f(b_1) + f(b_2) -f(\alpha) )/2\\
x_2 = (f(\alpha) + f(\gamma) -f(\beta))/2\\
x_3 = f(b_3)\\
x_4 = f(b_4)\\
x_5 = (f(\beta) -f(b_1) -f(b_3))/2\\
x_6 =(f(\beta) -f(b_2) -f(b_4))/2\endgather
$$

The solutions $x_i$'s are in $\bold Z_{\geq 0}$ by lemma 4.1, equations (4),
(6) and the assumption.

Subcase (2.3). $f(\beta) \geq f(b_1) + f(b_3)$ and $(f(\beta), f(b_2), f(b_4))
\in \Delta$.
$$\gather
x_1 + x_4 + x_5 = f(b_1)\\
x_1 + x_2 + x_3 + x_4 + 2x_6 = f(b_2)\\
x_3 = f(b_3)\\
x_2 + x_5 = f(b_4)\\
x_2 + x_3 + 2 x_4 + x_5 + 2x_6 = f(\alpha)\\
x_1 + x_3 + x_4 + x_5 + 2x_6 = f(\beta)\\
x_1 + x_2 + x_4 = f(\gamma)\endgather
$$
The solution is
$$\gather
x_1 = (f(b_1) + f(b_2) -f(\alpha) )/2\\
x_2 = (f(b_2) + f(b_4) -f(\beta))/2\\
x_3 = f(b_3)\\
x_4 = (f(\alpha) + f(b_1) -f(\beta) -f(b_4))/2\\
x_5 = (f(b_1) + f(b_4) -f(\gamma))/2\\
x_6 =(f(\beta) -f(b_1) -f(b_3))/2\endgather
$$
To show that the solutions are in $\bold Z_{\geq 0}$, it suffices to
show that $x_4 \in$ $\bold Z_{\geq 0}$ (the  rest of the $x_i \in \bold
Z_{\geq 0}$ follows from equations (4),(6), and the assumption). For $x_4$, we express $x_4$ as
$\frac{1}{2}((f(\alpha) + f(\beta) +f(\gamma)) -
(f(\beta) + f(b_2) +f(b_4))$. Thus $x_4$ is in $\bold Z_{\geq 0}$
by equations (4) and (6).

Subcase (3.3). Both $(f(\beta), f(b_1), f(b_3))$ and
$(f(\beta), f(b_2), f(b_4))$ are in $\Delta$.

The equation is,
$$\gather
x_1 + x_2 + x_3 + x_5 = f(b_1)\\
x_1 + x_4 + x_5 + x_6 = f(b_2)\\
x_3 + x_6 = f(b_3)\\
x_2 + x_4 = f(b_4)\\
x_2 + x_3 + x_4 + 2x_5 + x_6 = f(\alpha)\\
x_1 + x_2 + x_5 + x_6 = f(\beta)\\
x_1 + x_3 + x_4 + x_ 5 = f(\gamma)\endgather
$$
The solution is,
$$\gather
x_1 = (f(b_1) + f(b_2) -f(\alpha))/2\\
x_2 = ( f(b_4) + f(\beta) -f(b_2))/2\\
x_3 = (f(b_1) + f(b_3) - f(\beta))/2\\
x_4 = (f(b_2) + f(b_4) -f(\beta))/2\\
x_5 = (f( \alpha) - f(b_3) -f(b_4))/2\\
x_6 = (f(b_3) + f(\beta) -f(b_1))/2\endgather
$$
By equations (4), (6), the solutions are in $\bold Z_{\geq 0}$.

This ends the proof of the proposition for $I_{\delta}$. The proof of the
characterization of geometric functions on $\Cal S(\Sigma_{0,4})$ 
is the same as in \S3.
Indeed, first of all, the rational solutions of
$\Sigma_{i=1}^3 x_i = \max_{1 \leq i \leq 3;  1 \leq j \leq 4}
(2x_i,$$ 2y_j, \Sigma_{j=1}^4 y_j, x_1 + y_1 + y_2,$$ x_1 + y_3 + y_4,
x_2 + y_1 + y_3,$$ x_2 + y_2 + y_4, x_3 + y_1 + y_4,$$ x_3+ y_2 + y_3)$
 are dense in the solutions over  $\bold R_{\geq 0}$. Also
if we consider  $f(\alpha_1 \alpha_2)$ as an unknown in equation (4), it becomes
$x + a =\max(2x, x+b, c)$ where $c = f(\alpha_1 \alpha_2) + f(\alpha_2 \alpha_1)$ (by equation (5)). Thus, by lemma 3.1, we see that the corresponding 
lemma 3.2 holds for $\Sigma_{0,4}$. This shows 
that equations (4),(5) characterize the geometric functions.$\square$

\it Remark 4.1. \rm  The proof actually shows that except for at most four 
adjacent ideal triangles, equations (1), (4) become triangular
equalities $\sum_{i=1}^3 f(\alpha_i)$$ =$$ \max_{i=1}^3$\newline 
$(f(\alpha_i))$ when $f =I_{\delta}$ for $\delta \in \Cal CS$.

As a consequence of the discussion in the last paragraph and lemma 3.1(b),
we obtain,

{\bf Corollary 4.1}. \it 
(a) Suppose $\alpha_1 \perp_0 \alpha_2$ in $\Cal S(\Sigma_{0,4})$ so that
$(\alpha_1 \alpha_2, b_1, b_2)$ bounds a $\Sigma_{0,3}$, then
$f(\alpha_1) + f(\alpha_2) = \max(f(\alpha_1 \alpha_2), f(\alpha_2 \alpha_1),
f(b_1) + f(b_2), f(b_3) + f(b_4))$.

(b) Suppose $\alpha_1 \perp \alpha_2$ in $\Cal S(\Sigma)$. Then
$f(\alpha_1) + f(\alpha_2) = \max(f(\alpha_1 \alpha_2), f(\alpha_2 \alpha_2))$.

\rm

Combining propositions 3.1, 4.1, we obtain the following useful consequence.

{\bf Corollary 4.2.} \it Suppose $f : \Cal S(\Sigma) \to \bold R_{\geq 0}$
satisfies equations (1),(2),(4),(5) and $\alpha \perp \beta$, or $\alpha
\perp_0 \beta$  in $\Cal S(\Sigma)$. Then $f(\alpha^n \beta)$ is convex in $n \in
\bold Z$. Furthermore, there is an integer $N$ so that for $n \geq N$,
$f(\alpha^n \beta) = f(\alpha^{n-1}\beta) + f(\alpha)$
and 
$f(\beta \alpha^n) = f(\beta \alpha^{n-1}) + f(\alpha)$.

\rm

\it Remark 4.2. \rm   It is shown in \S8 that  $f(\alpha^n \beta)$ 
is convex in $n \in \bold Z$ for all $\alpha, \beta \in \Cal CS_0(\Sigma)$.
This seems to be an analogy with the fact that the geodesic length
functions are convex along the Thurston's earthquake paths ([Ker2], [Wo]).
I would like to thank P. Schmutz for  drawing my attention  to 
the convexity property. The operation $\alpha^n \beta$ is 
similar to the extension of the earthquake from the Teichm\"uller
space to the measured lamination space. See [Bo3], [Pa1], [Pa2] also
\S8 for more discussion.
 
Proof. Since $\alpha, \beta$ lie in an incompressible subsurface homeomorphic
to either $\Sigma_{1,1}$ or $\Sigma_{0,4}$, we may assume that $\Sigma$ $\cong$
 $\Sigma_{1,1}$ or $\Sigma_{0,4}$. We shall consider the case $\alpha
\perp_0 \beta$ only (the other case is similar and simpler). Let
$x_n = f(\alpha^n \beta)$, $n \in \bold Z$. Since $\alpha^n \beta \perp_0
\alpha$ with  $\alpha(\alpha^n \beta) = \alpha^{n+1} \beta$, we obtain 
following two equations for the sequence \{$x_n\}$ by equations (4),(5):
$$x_{n+1} + x_n + f(\alpha) =\max(2x_{n+1}, 2x_n, x_{n+1}+ b_{n+1}, x_n + b_n,
c) \tag 7 $$
where $b_{2n} = b_0$ and $b_{2n+1} =b_1$, and
$$x_{n+1} + x_{n-1}=\max(2x_n,  x_n + b_n, c). \tag 8$$
Now by $(8)$, $x_{n+1} + x_{n-1} \geq 2 x_n$. Thus  $f(\alpha^n \beta)$
is convex in $n$. To show that $x_n$ is linear in $n$ for $|n|$ large,
we shall consider $n >0$ only (the other case is similar).
By convexity,  $x_n$ is monotonic for $n$ large. If $\lim_n x_n = \infty$,
then $x_{n+1} \geq  x_n > \max(b_n, c, c/2)$ for $n$ large. Thus for $n$ large,
$(7)$ becomes, $x_{n+1} = x_n + f(\alpha)$. If $\lim_n x_n = L$ is a finite
number, take the limit
to the equations (7) and (8). We obtain:
$$ 2L + f(\alpha) = \max(2L, L + b_{\infty}, c)$$ and
$$ 2L = \max(2L, L + b_{\infty}, c)$$ where $b_{\infty} =\max(b_0, b_1)$.
Thus $f(\alpha) =0$. By  $(7)$, this shows $x_n = x_{n+1}$ for all $n$,
i.e., $f(\alpha^n \beta) = f(\alpha^{n-1}\beta)$ + $f(\alpha)$.
$\square$

\S 5. A  Reduction Proposition

The necessity of the conditions in theorem 1 is evident. To show the sufficiency,
we use induction on the norm $|\Sigma_{g,r}| = 3g+ r$ of a
surface $\Sigma_{g,r}$. By propositions 3.1 and 4.1,  
theorem 1 holds for $|\Sigma| =4$. 
If  $|\Sigma| \geq 5$, we decompose $\Sigma$ as a union of two incompressible
subsurfaces $X$, $Y$ so that $X \cap Y \cong \Sigma_{0,3}$ and $|X|, |Y|
< |\Sigma|$. For instance, if $g=0$, we take $X = \Sigma_{0,4}$, 
$Y=\Sigma_{0,r-1}$; if $g \geq 1$, we take $X = \Sigma_{1,1}$ and
$Y = \Sigma_{g-1, r+2}$. Note
that $\partial X \cap int(\Sigma)$ consists of a simple loop. See figure 9.

\midspace{0.1cm}
\centerline{\epsfbox{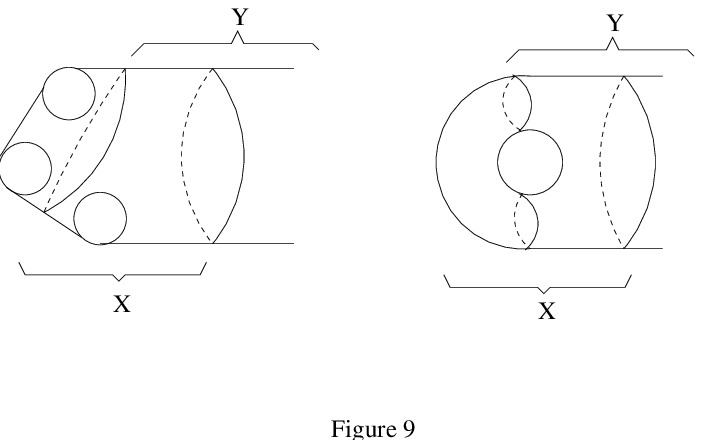}}
\midspace{0.1cm}

If $f:$$\Cal S(\Sigma)$ $\to \bold R$ satisfies equations (1), (2), (4), (5), then
$f|_{\Cal S(X)}$ and $f|_{\Cal S(Y)}$ again satisfy the same equations. By the
induction hypothesis, $f|_{\Cal S(X)} = I_{m_1}$ and $f|_{\Cal S(Y)} =I_{m_2}$
where $m_1 \in \Cal ML(X)$ and $m_2 \in \Cal ML(Y)$. Furthermore, by the
gluing lemma 2.1, there is $m \in \Cal ML(\Sigma)$ so that $m|_X = m_1$ and $m|_Y =m_2$.
Thus for $ h =I_m$, we have
$$ f|_{\Cal S(X) \cup \Cal S(Y)} = h|_{\Cal S(X) \cup \Cal S(Y)}  \tag 9$$
The goal of this and the rest of the sections \S6, \S7 is to show that
$f=h$ follows from (9).

{\bf Proposition 5.1.} \it Suppose $\delta \in \Cal S'(\Sigma)$ and $f, h:$
$\Cal S(\Sigma)$ $\to \bold R_{\geq 0}$ satisfy the equations (1),(2),(4),(5).
If $f(\alpha) = h(\alpha) $ for all $\alpha \in$ $\Cal S(\Sigma)$ with
$I(\alpha, \delta) \leq 2$, then $f =h$. \rm

Proof. We shall prove  that $f(\alpha) = h(\alpha)$ for $\alpha \in
\Cal S(\Sigma)$ by induction on the complexity
$(|\Sigma|, $$I(\alpha,$$ \delta))$ in the lexicographic order. By propositions 
3.1 and 4.1, it holds for $|\Sigma| \leq 4$. Assume now 
that $|\Sigma| \geq 5$ and
$\alpha \in \Cal S(\Sigma)$ so that $I(\alpha, \delta) \geq 3$. Take
$a \in \alpha$ and $d \in \delta$ so that $|a \cap d| = I(a,d)$. 
Fix an orientation on $a$. There are three
cases to be considered: (i) there are three intersection points $P_1$, 
$P_2, P_3$ in $a \cap d$ so that $P_1, P_2, P_3$ are adjacent along $d$ 
and their intersection signs are $(+,-,+)$ or $(-,+,-)$; (ii) there are three
adjacent (along $d$) intersection points $P_1, P_2, P_3$ in $a \cap d$ which
have the same intersection signs; and (iii) there are four adjacent 
intersection points $P_1, P_2, P_3, P_4$ in $a \cap d$ (along $d$) so 
that their intersection
signs are $(+,-,-,+)$ or $(-,+,+,-)$. See figure 10.

\midspace{0.1cm}
\centerline{\epsfbox{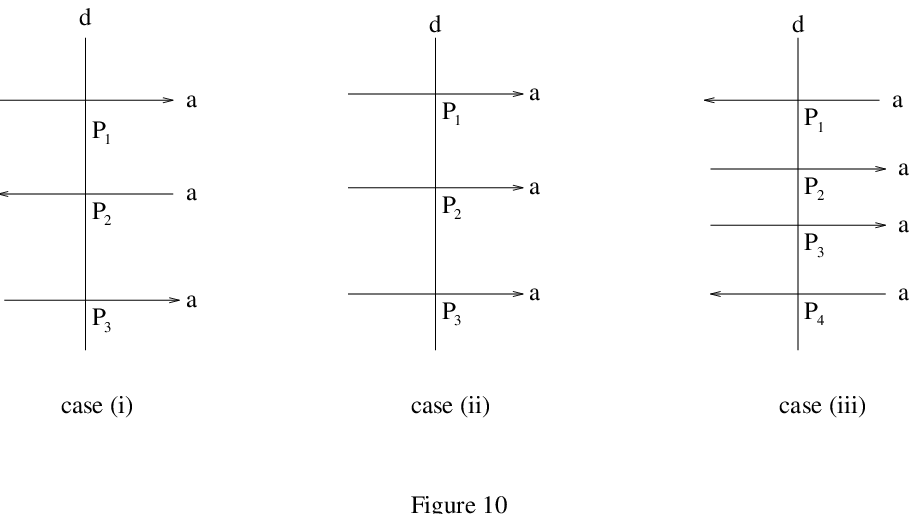}}
\midspace{0.1cm}

Case (i). There are two configurations  of $a \cup P_1P_3$ where $P_1 P_3$
is the arc in $d$ with end points $P_1$, $P_3$ so that $P_2 \in P_1 P_3$.
See figure 11(a), (b). These two cases are symmetric. Let us consider the
case  figure 11(b) only.

\midspace{0.1cm}
\centerline{\epsfbox{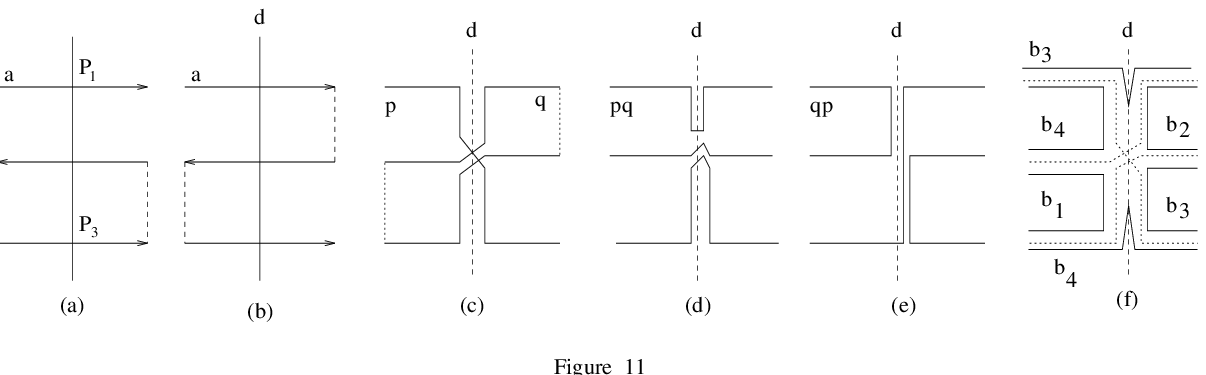}}
\midspace{0.1cm}

Let $p$, $q$ be two simple loops as indicated in figure 11(c). We have
$|p \cap q| =2$ and $p,q$ have zero algebraic intersection number. Since
$|a \cap d| =I(a, d)$, $p \perp_0 q$ (one way to see this is to
show that each component $b_i$ of $\partial N(p \cup q)$ is essential. Now
each $b_i$ is isotopic to a loop made by an arc $P_jP_{j+1}$, $j=1,2$, and
an arc along $a$ with end points $P_j$, $P_{j+1}$. 
Thus $b_i$ is essential since $|a \cap d| =I(a, d)$). 
Furthermore, $pq \cong a$, $I(qp, d),
I(p, d), I(q, d) < I(a,d)$ and $I(b_i, d) < I(a, d)$ where $\partial N(p \cup q)
= b_1 \cup b_2 \cup b_3 \cup b_4$ as shown in figure 11. Since $[p], [q],$
$[qp]$, and $ [b_i]$ are in $\Cal S(\Sigma)$ which have fewer intersections with $\delta$, by
the induction hypothesis, $f$ and $h$  have the same values on these seven
elements. By equation (5) applied to $[p] \perp_0 [q]$, we 
obtain $f(\alpha) = h(\alpha)$.

\it Remark 5.1. \rm In this case we conclude a stronger result that
if $d'$ is a curve system disjoint from $d$, then $I(p, d'), I(q, d'),
I(qp, d')$ and $I(b_i, d')$ $\leq I(a, d')$.

Case (ii). There are two configurations of $a \cup P_1P_3$ which are symmetric
(see figure 12(a), (b)).

\midspace{0.1cm}
\centerline{\epsfbox{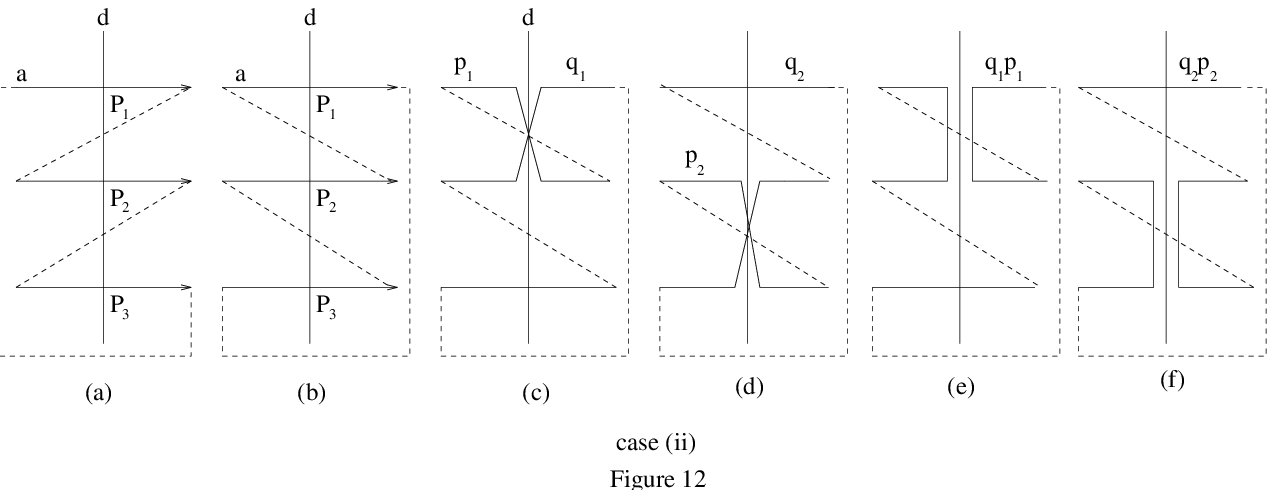}}
\midspace{0.1cm}

We shall consider the case in figure 12(b) only. Let $p_i$, $q_i$, $i=1,2$, be 
four simple loops as indicated in figure 12(c), (d). We have $p_i \perp q_i$,
$p_i q_i \cong a$, $|p_i \cap d|, |q_i \cap d|$ and $|q_i p_i \cap d| < I(a,d)$. One of the two curves $p_i$, say $p_1$, satisfies $|p_1 \cap d| < \frac{1}{2}
I(a,d)$ due to $|p_1 \cap d| + |p_2 \cap d| < I(a,d)$. Let
$b = \partial N(p_1 \cup q_1)$.  By the induction hypothesis and
 equation (2) for $p_1 \perp q_1$ 
that $f(\alpha) + f(q_1p_1) = \max(2f(p_1), 2f(q_1), f(b))$, 
$f(\alpha) =h(\alpha)$ follows from $f(b) = h(b)$.
The goal now is to show that $f(b) =h(b)$. Isotopy $p_1$ so that
$|p_1 \cap d| = I(p_1, d)$ and let $\Sigma'$ be $\Sigma -int(N(p_1))$. Then
the subsurface $\Sigma'$ is connected and incompressible since $p_1$ is
non-separating and essential. Furthermore, $[b] \in \Cal S(\Sigma')$ and
 $|\Sigma'| < |\Sigma|$. Thus
by the induction hypothesis, $f|_{\Cal S(\Sigma')} =I_{m_1}$ and 
$h|_{\Cal S(\Sigma')} = I_{m_2}$ for $m_1, m_2 \in
\Cal ML(\Sigma')$. We shall prove that $m_1 =m_2$. Thus in particular
$f(b) = h(b)$. To achieve this, let $d' = d \cap \Sigma'$ which is
a curve system consisting of $k$ arcs where $k =I(p_1, d) < \frac{1}{2} I(a,d)$.
By the induction hypothesis that  $f(\beta) = h(\beta)$ for $\beta$ with
$I(\beta, d) < I(a,d)$, we have $I_{m_1}( \beta) = I_{m_2}(\beta)$
for all $\beta \in \Cal S(\Sigma')$ so that $I(\beta, d') \leq 2k$ (
$ < I(a,d)$). Now $m_1 = m_2$ follows from the 
lemma below.

{\bf Lemma 5.1.} \it Suppose $F$ is a compact surface of negative Euler number
and $d$ is a curve system consisting of $k$ arcs. If $m_1$, $m_2 \in \Cal  ML(F)$
satisfy $I(m_1, \beta) = I(m_2, \beta)$ for all $\beta \in \Cal S(F)$
so that $I(\beta, d) \leq 2k$. Then $m_1 = m_2$. \rm

Proof. We use induction on $|F|$. If $|F| =4$, i.e., $F \cong \Sigma_{1,1}$
or $\Sigma_{0,4}$, then the result follows from propositions 3.1, and 4.1. 
Indeed,
each component $b \subset \partial F$ satisfies $I(b, d) \leq 2k$. Also
the ideal triangle $(\alpha, \beta, \gamma)$ in $\Cal S'(F)$ so that
$I_d(\alpha) + I_d(\beta) + I_d(\gamma)$ is minimal (among all such
triples) satisfies $I_d(\alpha), I_d(\beta), I_d(\gamma) \leq 2k$ (by the
proof of propositions 3.1, 4.1). Thus
$m_1 =m_2$ in this case.

If $|F| \geq 5$, we construct an ideal triangulation $ [t_1 ...  t_n]$ of 
$F$ so that $d \cong t_1^{k_1} ... t_n^{k_n}$ where $k_i \in
\bold Z_{\geq 0}$ and $\Sigma_{i=1}^n k_i = k$. There
are two  components of t, say $t_1 $and  $t_2$, so that  each of them is
non-separating. Indeed, it is known that any Fenchel-Nielsen system
on a surface $\Sigma_{g,r}$ must contain at least $g$ many non-separating
simple loops. By doubling  the surface $F$ and the ideal triangulation $t$,
we obtain the two non-separating arcs above. Let $X_i = F -int(N(t_i))$, $i=1,2$.
By the choice of $t_i$, each $X_i$ is connected and incompressible in $F$. Since
$|F| \geq 5$, the Euler number $\chi(X_i)$  of $X_i$ is negative. 
 Furthermore,
$\chi (X_1 \cap X_2) = \chi(F -int(N(t_1 \cup t_2))) < 0$ if $F \neq \Sigma_{1,2}$.
Now consider the restrictions $I_{m_j}|_{\Cal S(X_i)}$, $i,j=1,2$ and
$d_i = d|_{X_i}$. The curve system $d_i$ consists of at most $k$ arcs
and $I_{m_1}(\beta) = I_{m_2}(\beta)$ for all $\beta \in \Cal S(X_i)$
 with $ I(\beta,  d_i) \leq 2k$ by the hypothesis. Thus, by the induction hypothesis
applied to $X_i$ with respect to $d_i$, we have $m_1 |_{X_i} = m_2|_{X_i}$
for $i=1,2$.  

Now if $F \neq \Sigma_{1,2}$, then $ X_1 \cap X_2$ contains an incompressible
$\Sigma_{0,3}$. By the gluing lemma 2.1, we obtain $m_1 =m_2$.

If $F = \Sigma_{1,2}$ with $\partial F = b_1 \cup b_2$, then $X_1 \cap X_2
\cong \Sigma_{0,2}$ as in figure 13 (there are  four cases depending on the
locations of $t_i$ with respect to $b_j$'s).

\midspace{0.1cm}
\centerline{\epsfbox{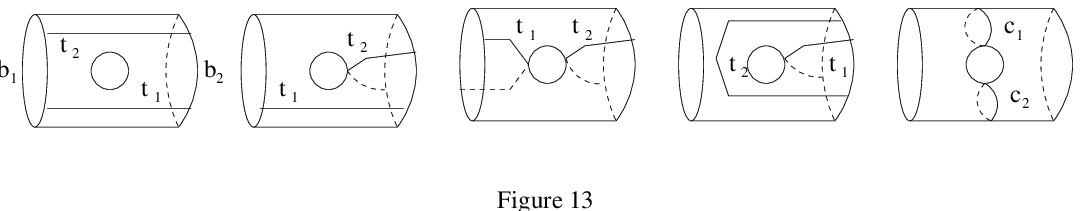}}
\midspace{0.1cm}

Let $[c_i] \in \Cal S'(X_1) \cup \Cal S'(X_2) $ so that $\{ c_1, c_2\}$
forms a Fenchel-Nielsen system. Then
each $m \in \Cal ML(F)$ is determined by its intersection numbers with
$b_1$, $b_2$, $c_1$, $c_2$ and the twisting coordinates at $c_1$, $c_2$
(the Dehn-Thurston coordinates). Now $m_1|_{X_i} = m_2|_{X_i}$ shows
that their twisting coordinates are the same at $c_1, c_2$. Thus $m_1 =m_2$. $\square$.

Case (iii).
There are four adjacent intersection points $P_1, P_2, P_3, P_4$ in $a \cap d$
so that their intersection signs are either $(+,-,-,+)$ or $(-,+,+,-)$.
Let $P_1P_4$ be the arc in $d$ with end points $P_1, P_4$ so that $P_2 
\in P_1P_4$ and let $P_iP_j$ be the arc in $d$ with ends $P_i, P_j$ so that $P_iP_j \subset P_1P_4$.
There are six possible configurations of $a \cup P_1P_4$ as  shown in
figure 14. Due to symmetry, it suffices to consider the cases (3.1), (3.2),
(3.3) in figure 14.

\midspace{0.1cm}
\centerline{\epsfbox{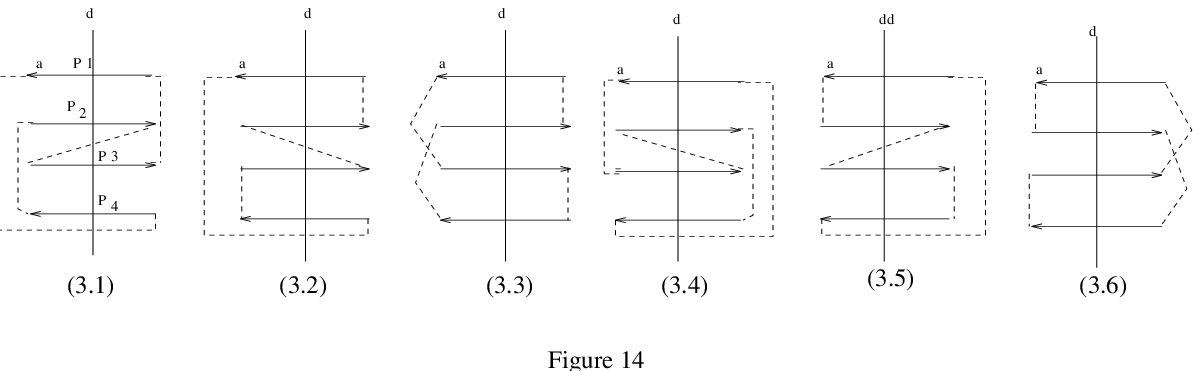}}
\midspace{0.1cm}

Note that $P_1P_2$ (resp. $P_3P_4$) approaches its end points from the
same side of $a$ and $P_1P_2$, $P_3P_4$ approach $P_1, P_3$ from 
different sides of $a$. Thus, $\Sigma'$ =$N(a \cup P_1P_2 \cup P_3P_4) \cong \Sigma_{0,4}$.
Furthermore, since $|a \cap d| = I(a, d)$, $\Sigma'$ is incompressible in $\Sigma$ so 
that $\partial \Sigma' = b_1 \cup b_2 \cup b_3 \cup b_4$ satisfies
$|b_i \cap d| < |a \cap d|$, $i=1,2,3,4$. In particular, $f(b_i) = h(b_i)$,
$i=1,2,3,4.$

Subcase 3.1. Let $\beta, \gamma \in \Cal S'(\Sigma)$ be as shown in figure
15 with $\beta \perp_0 \gamma$.

\midspace{0.1cm}
\centerline{\epsfbox{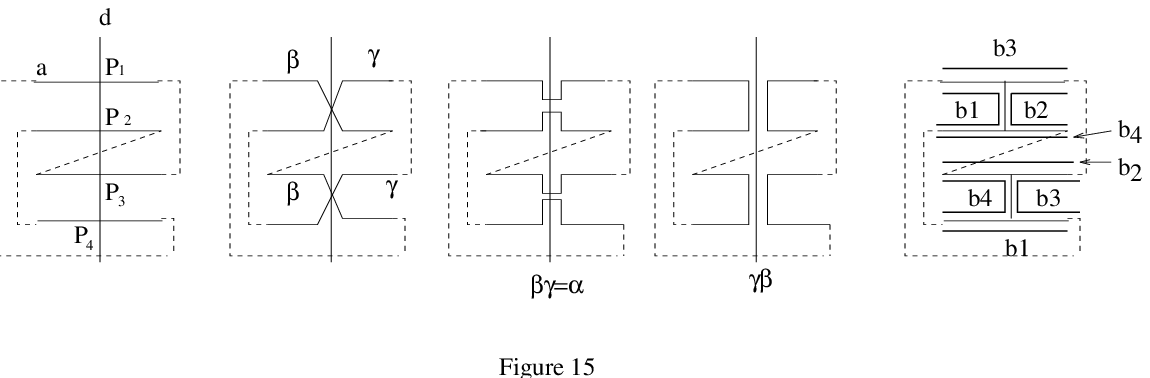}}
\midspace{0.1cm}

We have $\alpha =\beta \gamma$, $I(\beta, d), I(\gamma, d), I(\gamma \beta, d)
< I(a,d)$. Thus, by the induction hypothesis and equation (5) for
$\beta \perp_0 \gamma$, $f(\alpha) = h(\alpha)$.

Subcase 3.2. Let $\beta \in \Cal S'(\Sigma)$ be as shown in figure 16 with $\beta \perp_0 \alpha$
and $I(\beta, d) < I(\alpha, d)$.

\midspace{0.1cm}
\centerline{\epsfbox{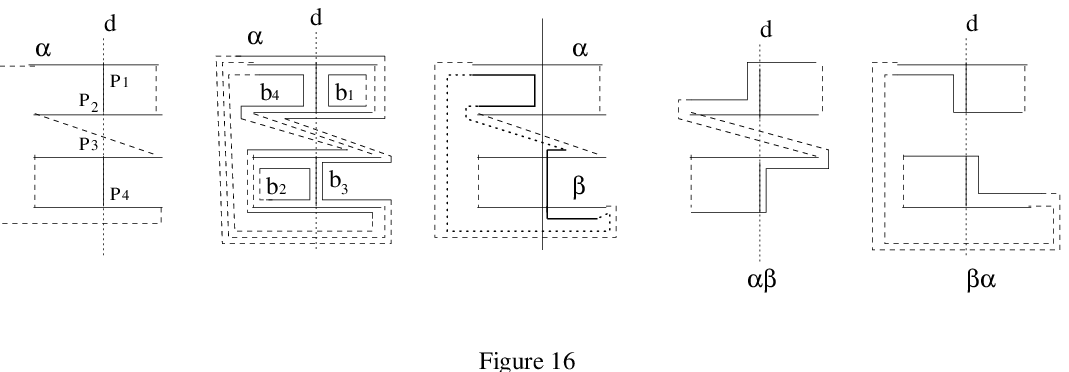}}
\midspace{0.1cm}

To simplify notations, in the rest of this section, 
we shall use $\alpha \beta$, $\beta \alpha$ to denote the
simple loops representing  them as shown in the figures.

The curves $\alpha \beta$ and $\beta \alpha$ are as indicated in figure 16. By corollary 4.1 applied to $\alpha_1 = \alpha$ and $\alpha_2 =\beta$, 
we see that $f(\alpha) = h(\alpha)$ follows from the claim below.

Claim.  $f(\alpha \beta) = h(\alpha \beta)$ and
$f(\beta \alpha) = h(\beta \alpha)$.

Proof of the claim. To show $f(\alpha \beta) = h(\alpha \beta)$, we observe that
the three adjacent intersection points $P_1, P_2, P_3$ in $\alpha \beta \cap d$
(along $d$) have alternating intersection signs. Thus the reduction process of
case (i) applies.

\midspace{0.1cm}
\centerline{\epsfbox{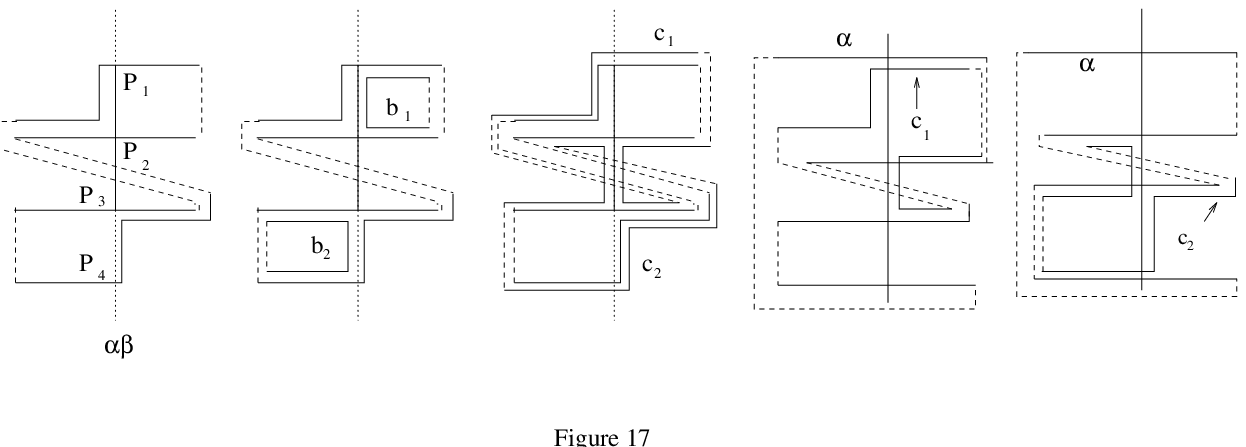}}
\midspace{0.1cm}

Consider the subsurface $X = N(\alpha \beta \cup P_1P_3) \cong \Sigma_{0,4}$.
The boundary components of $X$ are isotopic to $b_1, b_2, c_1, c_2$ where
$b_1, b_2 \subset \partial \Sigma' =\partial N(a \cup P_1P_2 \cup P_3P_4)$
and $c_i \perp
a$. Thus $c_i$'s are essential simple loops in $\Sigma$ and $X$ is incompressible in $\Sigma$.
Furthermore,  $I(c_i, d) < I(a,d)$ as shown in figure 17. Consider two classes $\beta', \gamma'
\in \Cal S(X)$ with $\beta' \perp_0 \gamma'$ as in figure 18. 
We have
$\beta' \gamma' = \alpha \beta$ and $I(\beta', d), I(\gamma', d), I(\gamma'
\beta', d) < I(a,d)$ as shown in figure 18. By the inductin hypothesis, $f, h$ have the same values on 
$\{\beta', \gamma', \gamma'\beta', b_1, b_2, c_1, c_2\}$. By equation (5), 
$f(\beta' \gamma') =h(\beta' \gamma')$, i.e.,  $f(\alpha \beta) = h(\alpha \beta)$.

\midspace{0.1cm}
\centerline{\epsfbox{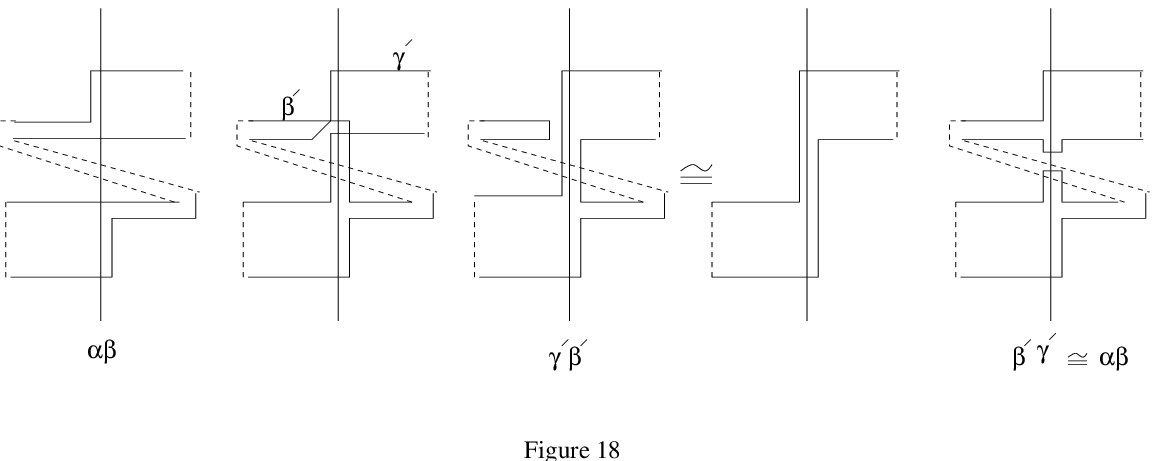}}
\midspace{0.1cm}

The proof of $f(\beta \alpha) = h(\beta \alpha)$ is similar. Take
$Y = N(\beta \alpha \cup P_1 P_3)$. Then $\partial Y$
consists of simple loops isotopic to $b_1, b_2$, $c_3, c_4$ as shown in 
figure 19.

\midspace{0.1cm}
\centerline{\epsfbox{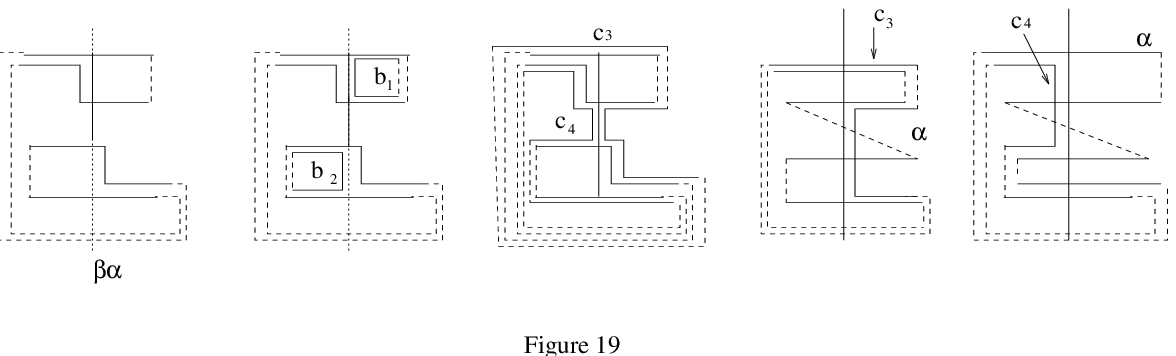}}
\midspace{0.1cm}

Note that $c_j \perp a$,  $I(c_j, d) < I(a,d)$ ($j=3,4$) 
and $b_i \subset \partial \Sigma'$. 
Thus $Y$ is incompressible in $\Sigma$. Now consider
$\beta', \gamma'$ $\in \Cal S(Y)$ as in figure 20. Then $\beta' \perp_0
\gamma'$, $\beta' \gamma' = \beta \alpha$, and $I(\beta', d), I(\gamma', d),
I(\gamma' \beta', d) < I(a,d)$. Thus by the induction hypothesis
and equation (5), $f(\beta \alpha) = h(\beta \alpha)$.

\midspace{0.1cm}
\centerline{\epsfbox{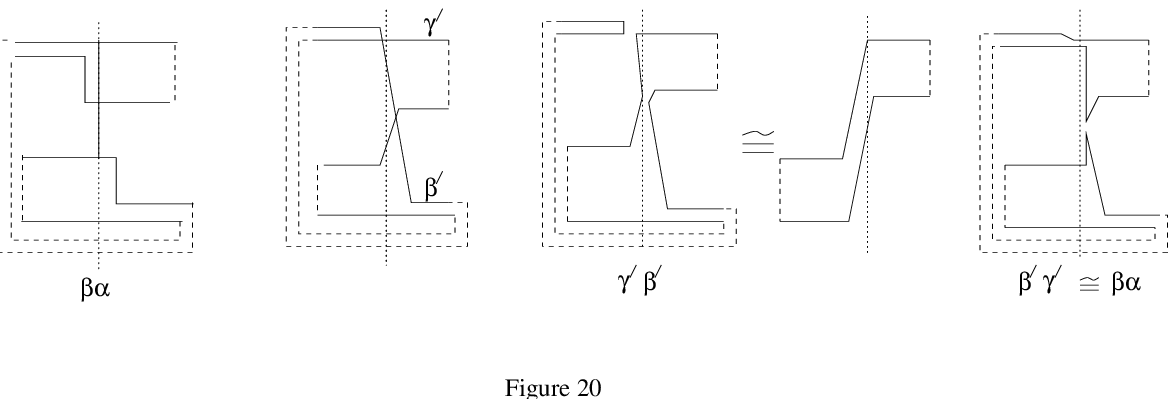}}
\midspace{0.1cm}

Subcase (3.3). This case is similar to the subcase (3.2). Let $\beta \in$ $\Cal S'(\Sigma)$
be as shown in figure 21 where $\beta \perp_0 \alpha$ and $I(\beta, d) < I(a,d)$.

\midspace{0.1cm}
\centerline{\epsfbox{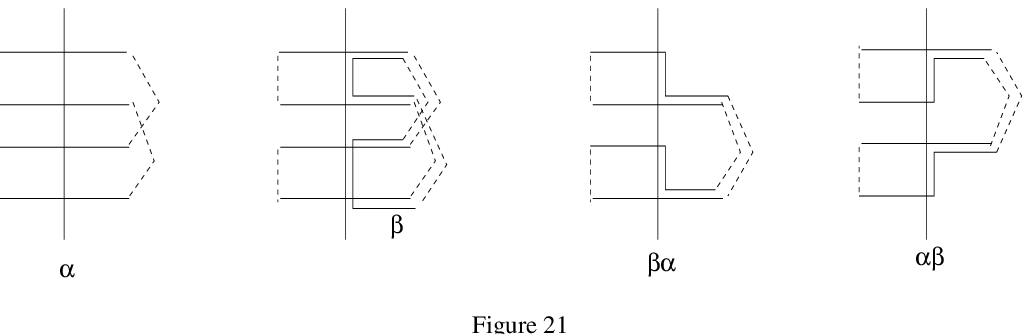}}
\midspace{0.1cm}

By the same argument as in subcase (3.2), it suffices to show $f(\alpha \beta)
=h(\alpha \beta)$ and $f(\beta \alpha) = h(\beta \alpha)$. We prove $f(\alpha
\beta) =h(\alpha \beta)$ below (the other case follows by symmetry).
Consider $X = N(\alpha \beta \cup P_1P_3)$. Then $\partial X \cong
b_1 \cup b_2 \cup c_1 \cup c_2$ as shown in figure 22 so that each of 
the component
has fewer intersection points with $d$. Furthermore, and $c_i \perp a$ and
$b_i$ $\subset \partial N(a \cup P_1P_2 \cup P_3P_4)$. Thus 
$X$ is incompressible in $\Sigma$.

\midspace{0.1cm}
\centerline{\epsfbox{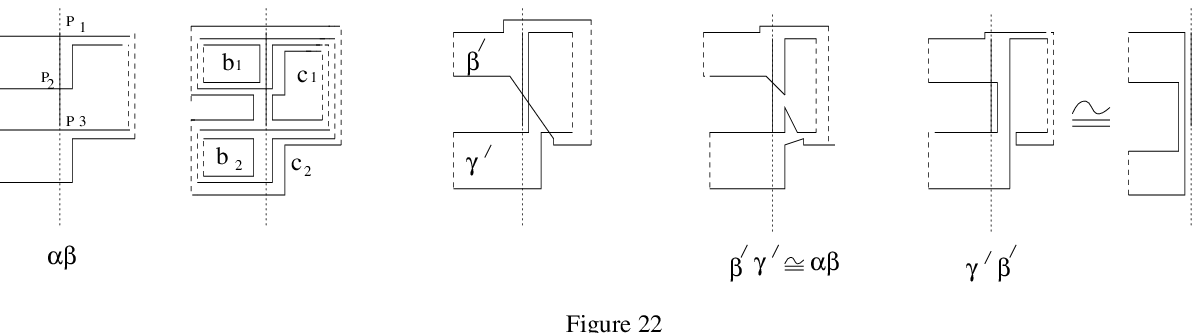}}
\midspace{0.1cm}

Now choose $\beta' \perp_0 \gamma'$ in $\Cal S(X)$ so that $\beta' \gamma'
= \alpha \beta$, and $I(\beta', d),$$I(\gamma', d),$$I(\gamma' \beta', d)$
$ < I(a,d)$ as shown in figure 22.
Thus by the induction hypothesis and equation (5), $f(\alpha \beta) = h(\alpha \beta)$.
$\square$

If $a, d$ are simple loops on $\Sigma_{0,r}$ so that $|a \cap d| \geq 3$, then
any three adjacent intersection points in $a \cap d$ have alternating
intersection signs. Thus, the cases (ii), (iii) in the proof of proposition 5.1
do not occur. Combining this observation and remark 5.1, we obtain,

{\bf Corollary 5.1.}  \it Suppose $\delta_1 ... \delta_k \in \Cal CS(\Sigma_{0,r})$ forms a Fenchel-Nielsen system on the surface and $f,g$ $:\Cal S(\Sigma_{0,r}
) \to \bold R_{\geq 0}$ satisfy equations (4),(5). If $f(\alpha) = h(\alpha)$
for all $\alpha \in \Cal S(\Sigma)$ with $I(\alpha, \delta_i) \leq 2$
for all $i$, then $f = h$. \rm

For surface $\Sigma_{1,2}$, the situation is more complicated.

{\bf Corollary 5.2.} \it Suppose $[a][b] \in \Cal S'(\Sigma_{1,2})$
forms a Fenchel-Nielsen system on the surface so that $a$ is separating.
If $f,h:\Cal S \to \bold R$  satisfy equations (1),(2),(4),(5) and
$f(\alpha) = h(\alpha)$ for $\alpha \in \{ \alpha \in \Cal S:
$ either $I(\alpha, a)I(\alpha, b) =0$ or $\alpha \perp_0 [a]$ and $\alpha
\perp [b]$\}, then $f=h$. \rm

Proof. By proposition 5.1 applied to $\delta = [a]$, it suffices to show that
$f(\alpha) = h(\alpha)$ for $\alpha \perp_0 [a]$. Let $k = I(\alpha, b)$.
We shall prove the proposition by induction on $k$. The assumption shows
that $f(\alpha) = h(\alpha)$ for $k \leq 1$. If $k \geq 2$, let $X$ 
be the subsurface $\Sigma_{1,1}$ bounded by $a$ and let $x \in \alpha$
so that $|x \cap b| = k, |x \cap a| = 2$. Thus $x \cap X$ consists
of an arc. Let $s_1$ be an essential simple loop in $X$ so that $x  \cap
s_1 = \emptyset$ and $s_2$ be an essential simple loop  in $X$
so that $s_2 \perp s_1$
and $s_2 \perp x$.

\midspace{0.1cm}
\centerline{\epsfbox{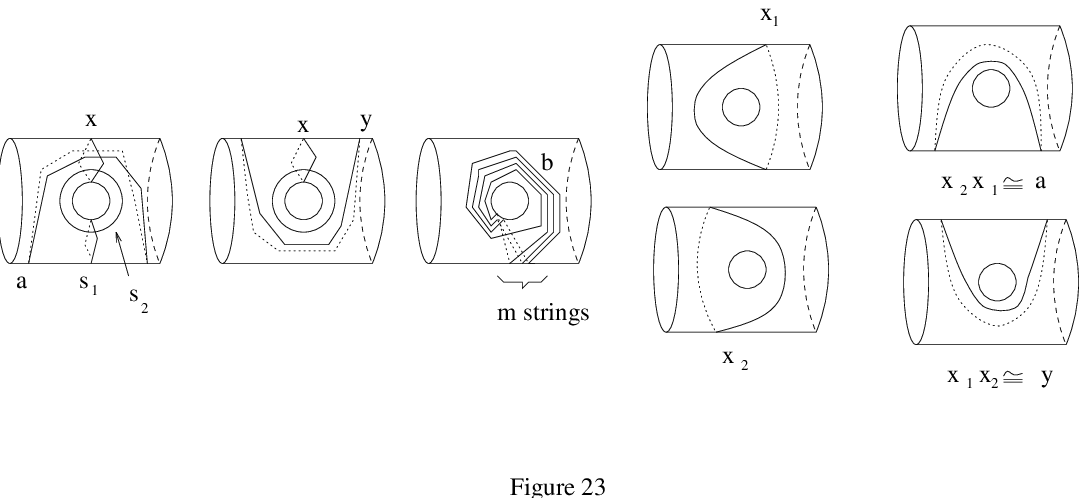}}
\midspace{0.1cm}

Since $b \subset X$ and $|b \cap x| = I(b,x) = k$, we may
express $b$ as $s_2^k s_1^m$ where $m \in \bold Z$. The simple loop $s_2$ is
not unique up to isotopy since we may replace it by $s_1^n s_2$. Replacing 
$s_2$ by an appropriate $s_1^n s_2$, we may write  $b \cong s_2^k s_1^m$
where $|m| \leq k/2$.  Thus $I(s_2, b) = |m| \leq k/2$. Let $y = \partial
N(x \cup s_2)$. Then $I(y, b) \leq 2|m| \leq k$. Assume for definiteness that
$m >0$. Then for $i=1,2$, $I(s_2^i x, b) = I(s_2^i x, 
(s_2^i s_1)^{k_1}(s_2^{i-1}s_1)^{m_1}) < k$. To see this, we note that $I(s_2^i x, s_2^i s_1) = 0$, 
$I(s_2^i x, s_2^{i-1} s_1) = 1$ and $I(s_2^is_1, s_2^{i-1}s_1) = 1$.
Now we express $b$ as $(s_2^i s_1)^{k_1}(s_2^{i-1}s_1)^{m_1}$ with
$ m_1 = k -im $ and $|m_1| < k$. Thus the result follows.  On the  other hand 
$s_2 \perp s_2 x$, and $x \cong (s_2x)s_2$. By the induction hypothesis, 
$f$, $h$ have the same values on $s_2x$, $s_2$, $s_2(s_2x) = s_2^2x$. 
Thus $f(x) = h(x)$ follows
from $f(y) = h(y)$ by equation (2) ($y \cong \partial N(s_2 \cup s_2 x)$). To show
that $f(y) = h(y)$, we consider $x_1  = s_1 s_2 x$ and $x_2 = xs_2 s_1$ 
as shown in figure 23. We have $x_1 \perp_0 x_2$, $x_1 x_2 \cong y$, 
$x_2 x_1 \cong a$,
$ N(x_1 \cup x_2) \cong \Sigma - int(N(s_2))$. Let
$\partial \Sigma_{1,2} = b_1 \cup b_2$.  By the
construction, $x_i \perp_0 a$ and $|x_i \cap b| \leq k/2$, $i=1,2$. This
shows that $f$ and $h$ have the
same values at $x_1$, $x_2$, $x_2x_1$, $s_2$, $b_1$, $b_2$ by the
induction hypothesis. Thus $f(y) = h(y)$
by equation (5) for $x_1 \perp_0 x_2$. $\square$.

\S 6. The Two-holed Tours and the Five-holed Sphere

We prove theorem 1 for surfaces $\Sigma$ = $\Sigma_{1,2}$ and $\Sigma_{0,5}$ in this section.

Choose two disjoint essential simple loops $a,b$ in $\Sigma$ so that
(1) $a$ is separating and (2) \{$a,b\}$ forms a  Fenchel-Nielsen system on $\Sigma$
as in figure 24 (a), (b).

\midspace{0.1cm}
\centerline{\epsfbox{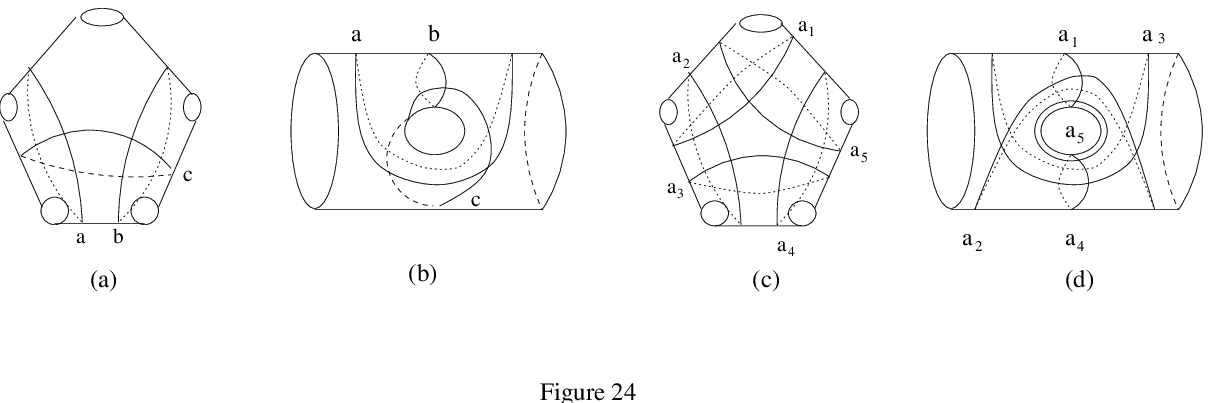}}
\midspace{0.1cm}

{\bf Lemma 6.1.} \it Let $c$ be an essential simple loop in $\Sigma$ so that
either $c \perp_0 a$, $c \perp_0 b$ or $c \perp_0 a$, $c \perp b$. 
Suppose $f,h : \Cal S(\Sigma) \to \bold R$ satisfy equations (1),
(2),(4),(5) and $f(\alpha) = h(\alpha)$ for
all $\alpha$ in $\{ \alpha \in \Cal S :$ either $I(\alpha, a)
$$ I(\alpha, b) =0$, or $\alpha$$ = [c],$$ [ac],$$ [cb], [acb]$\}. Then $f =h$. \rm

Proof. By corollaries  5.1 and 5.2, it suffices to show that $f(\alpha) = h(\alpha)$ for  either $\alpha \perp_0 [a]$,  $\alpha \perp [b]$ or $\alpha \perp_0
[a]$, $\alpha \perp_0 [b]$. By comparing the Dehn-Thurston coordinate at $\{a,b\}$, we have $\alpha =[a^ib^jc]$ for
some $i,j \in \bold Z$. Take $x \in \alpha$ so that $|x \cap a| =I(x,a)$
and $|x \cap b| = I(x,b)$. Since  either $x \perp_0 a$,  $x \perp b$ or
$x \perp_0 a$, $x \perp_0 b$, $\partial N(x \cup a)$ and $\partial N(x \cup b)$ are
either  disjoint from $a$ or from $b$. In particular, $f, h$ have the
same values at these boundary components. Thus to show $f(x) =h(x)$, by 
equation (5) applied to $x  \cong a (xa)$,
it suffices to show, for instance,  $f(xa) = h(xa)$, $f(xaa) = h(xaa)$.
We shall prove this by induction on $||\alpha|| = |i| + |j|$. If
$||\alpha|| =1$, then $x \cong ac, ca, bc, cb$. Now $f(ca) = h(ca)$ follows
from $f(ac) = h(ac)$, $f(a) = h(a)$, and $f(c) = h(c)$ (by equation (5)). Similarly,
we have $f(bc) = h(bc)$. If $||\alpha||=2$, then $x \cong abc, acb, bca, cab$.
To show for instance that $f(abc) = h(abc)$, we write $abc \cong  b(ac) $.
 Now $f,h$ have the same values on \{$ac, b, (ac)b)$\}. 
Thus by
equations (2) or (5), we have $f(abc) = h(abc)$. By the same argument
we see that $f$ and $h$ have the same values at $bca$, $cab$. 
Suppose now
that $||\alpha|| \geq 3$. Then one of the numbers $|i|$ or $|j|$ is at least
2. Say, $|i| \geq 2$. For definiteness, we assume that $i \geq 2$ (the other
case $i \leq -2$ is similar). Write  $x = a^i b^j c = a(a^{i-1} b c): = ay$
where $||[y]|| < ||\alpha||$. Furthermore, $ya = a^{i-2}b^j c$ has
norm $|| [ya]|| < ||\alpha||$. Thus $f,h$ have the same values at
$\{a, y, ya\}$ by the induction hypothesis. 
We obtain $f(x) = h(x)$ by equation (5) (or equation (2) in case
$|i| \leq 1$ and $|j| \geq 2$). $\square$

We now begin the  proof of theorem 1 for $\Sigma_{1,2}$ and $\Sigma_{0,5}$. Given a  non-zero function
$f: \Cal S \to \bold R_{\geq 0}$ satisfying equations (1),(2),(4),(5), we  choose a
pair of elements $[a], [b] \in \Cal S'(\Sigma)$ so that (1) $a$ is
separating, (2) $\{a,b\}$ forms a Fenchel-Nielsen system, and (3)
$f(a)f(b)$ is non-zero. To see that condition (3)
can be realized, we use the fact that if a geometric function
$k: \Cal S(F) \to \bold R$  takes non-zero values at $\partial F$
then $k|_{\Cal S'(F)} \neq 0$. 

By the reduction process in \S5, we construct a measured 
lamination $m \in \Cal ML(\Sigma)$ so that $ f(\alpha) = I_m(\alpha)$
for all $\alpha$ satisfy $I(\alpha, a) I(\alpha, b) = 0$. Call $h = I_m$
for simplicity. By lemma 6.1, it suffices to find $[c] \in \Cal S$ so that
$c \perp_0 a$, $c \perp b$ or $c \perp_0 b$ and $f, h$ have the
same values at $\{c, ac, cb, acb\}.$ 

We shall consider $\Sigma$ =$\Sigma_{1,2}$ and $\Sigma_{0,5}$ separately.

Case 1. $\Sigma$ = $\Sigma_{1,2}$.

{\bf Lemma 6.2.} \it  Suppose $[a'], [b'] \in \Cal S(\Sigma_{1,2})$
so that (1) $a' \perp b'$, $a' \perp_0 a$, $b' \perp b$,
$a' \cap b = a \cap b' = \emptyset$, and (2) $f(aa') + f(b'bb) 
< f(a') + f(b')$. Then  $f(a'b') = h(a'b') = f(a') + f(b')$ \rm

Proof. First by figure 25, we have 
$aa' \perp b'bb$ and $b'a' \cong (aa')(b'bb)$.

\midspace{0.1cm}
\centerline{\epsfbox{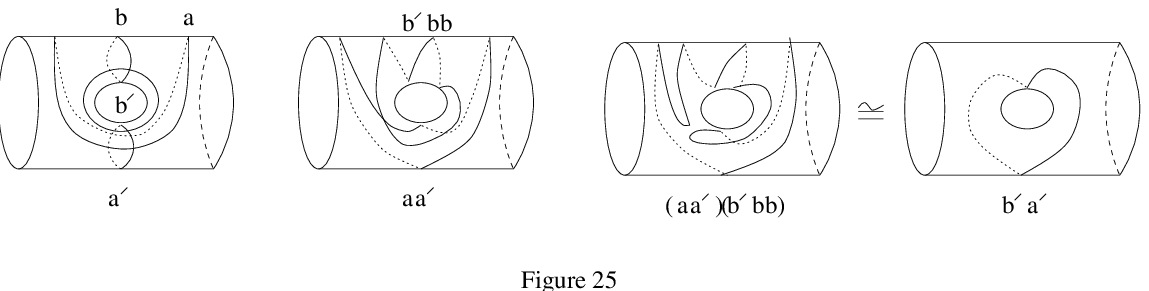}}
\midspace{0.1cm}

By the triangular inequality $f(\alpha \beta) \leq f(\alpha) + f(\beta)$
whenever $\alpha \perp \beta$ or $\alpha \perp_0 \beta$, we obtain $f(b'a') \leq f(aa') + f(b'bb)
< f(a') + f(b')$. By corollary 4.1(b), $f(a') + f(b') = \max(
f(a'b'), f(b'a'))$. Thus $f(a'b') = f(a') + f(b')$. Since $f,h$ have the
same values at simple loops disjoint either from $a$ or from $b$, the
same argument applies to $h$. We conclude that $h(a'b') = h(a') + h(b')$
$= f(a') + f(b') = f(a'b')$. $\square$

To prove theorem 1 for $\Sigma_{1,2}$, take $a_1 \perp_0 a$, $b_1 \perp b$,
$a_1 \perp b_1$ so that $a_1 \cap b = a \cap b_1 = \emptyset$ (as shown in
figure 26). For any integers $n$, $m$, $a' = a^n a_1$
and $b' = b_1 b^m$ satisfy the condition (1) in lemma 6.2.
By corollary 4.2, we may replace $a_1$ by $a_1 a^n$ and
$b_1$ by $b^n b_1$ for some large $n$ so that after the replacement,
$$ f(aa a_1) < f(a a_1) < f(a_1) \tag 10$$ and
$$ f(b_1 bbb) < f(b_1 bb) < f(b_1 b) < f(b_1) < f(b_1 b) \tag 11 $$

The same inequalities also hold for $h$ 
since $f$ and $h$ have the same values at the simple loops disjoint
either from $a$ or from $b$.

Take $c = a_1 b_1$.  Applying lemma 6.2  to $f,h$ with $a' = a_1$ and $
b'= b_1$, we obtain $f(c) = h(c)$ (the conditions in 
the lemma are satisfied due to (10) and (11)).

Now $cb \cong a_1 ( b_1 b)$. Take $a' = a_1$ and $b' = b_1 b$ in lemma 6.2.
We obtain $f(cb) = h(cb)$.
Also $ac \cong (a a_1) b_1$. Take $a' = aa_1$ and $b' = b_1$ in lemma 6.2.
We obtain $f(ac) = h(ac)$. 

Finally, note that $acb \cong (aa_1)(b_1b)$ as shown in figure 26.

\midspace{0.1cm}
\centerline{\epsfbox{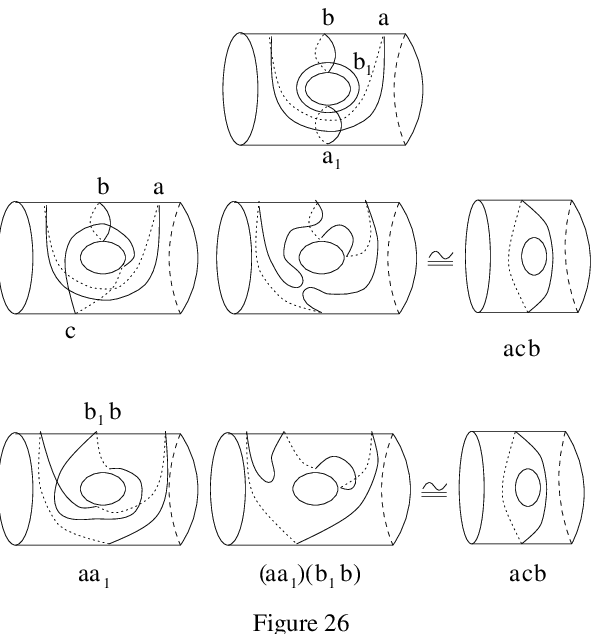}}
\midspace{0.1cm}

Take $a' = aa_1$ and $b' =b_1 b$ in lemma 6.2. We obtain $f(acb) = h(acb)$.

Case 2. $\Sigma$ = $\Sigma_{0,5}$. Suppose $\partial \Sigma_{0,5} =
\{ \partial_1,..., \partial_5\}$ and let $M 
=\max\{f(\partial_i): i=1,2,3,4, 5\}$.
First, we  make $\min(f(a), f(b))$ arbitrary large by choosing
different pairs $\{a,b\}$. To see this, choose
$a_1 \perp_0 a$, $a_1 \cap b =\emptyset$ and $f(a_1) \neq 0$. 
Now replace $a$ by $a a_1^n$
for a large n.
Then corollary 4.2 shows that $f(a a_1^n)$ growth linearly in $n$.
By repeating the replacement inside the surface $\Sigma  - int(N(a a_1^n))$,
we can make $f(b)$ large as well. Thus we may assume that
$\min(f(a), f(b))$ $> 4M$

{\bf Lemma 6.3.} \it Suppose $a', b'$ are two essential simple loops
in $\Sigma_{0,5}$ so that (1) $a' \perp_0 b'$,
$a \perp_0 a'$, $b \perp_0 b'$, $a' \cap b = a \cap
b' = \emptyset$, and (2) $f(a') + f(b') > f(aa') + f(b'b) + 2M$.
Then $f(a'b') = f(a') + f(b')$.\rm

Proof. Let $ a'' = aa',  b'' = b'b$, $x = a'b', x' = b'a'$,
$y = a''b''$, and $y' = b''a''$. Then
we have:  $I(y, a') = I(y, b')= 0$ and $I(x' , a'') = I(x' , b'')  =0$
by figure  27. This shows that
$N( a' \cup b')$ is isotopic the  subsurface $\Sigma_{0,4}$ bounded by $y$, and 
$N(a'' \cup b'')$ is isotopic to the subsurface
 $\Sigma_{0,4}$ bounded by $x'$. Now by corollary 4.1(a), we have
$$ f(a') + f(b')  = \max(f(x), f(x'), f(y) + M_1', 2M'_2)  \tag 12$$
and $$ f( a'') + f( b'') = \max(f(y), f(y'), f(x') + M_1'', 2M''_2) \tag 13$$
where $0 \leq M_i', M_i'' \leq M$, $i=1,2$. 
By (13), $f(y)  \leq f( a'') + f(b'')$. Thus
$f(y) + M_1' \leq  f(a'') + f(b'') + M_1' < f(a') + f(b')$
by the assumption.
Also by (13), $f(x') \leq  f(a'') + f(b'') < f(a') + f(b')$
by the assumption again. Finally  by the assumption,
$f(a') + f(b') > 2M_2'$.  Thus, (12) becomes $f(x) =  f(a') + f(b')$
$\square$

\midspace{0.1cm}
\centerline{\epsfbox{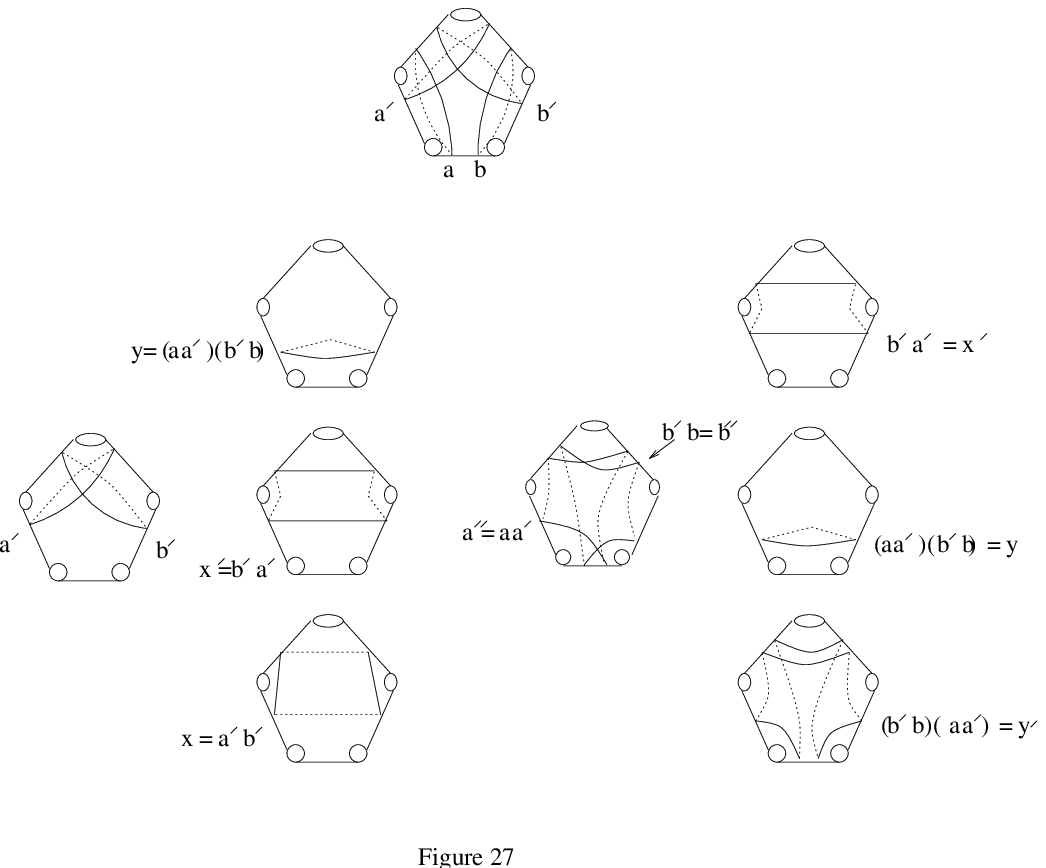}}
\midspace{0.1cm}

Now we apply the lemma to finish the proof of theorem 1. Take
$a_1 \perp_0 a$, $b_1 \perp_0 b$, $a_1 \perp_0 b_1$, and
$a \cap b_1 = a_1 \cap b = \emptyset$. 
Then for any $n$, $m$, $a' = a^n a_1$
and $b' = b_1 b^m$ satisfy the condition (1) in lemma 6.3. 
By corollary 4.2, 
replacing $a_1$ by $a_1a^n$ and $b_1$ by $b_nb_1$ for $n$ large, 
we  may assume that
$$ f(a^ia_1) > f(a^{i+1} a_1) + 2M,  \quad \quad f(b_1 b^i) > f(b_1 b^{i+1}) 
+ 2M, \quad i=1,2,3 \tag 14$$
Take $c =  a_1 b_1$. Then by lemma 6.3 applied to $a' = a_1$ and
$b' = b_1$, we obtain $f(c) = h(c)$.
Take $a' = aa_1$ and $b'=b_1$ in lemma 6.3. We obtain $f(ac) = h(ac)$.
Take $a' = a_1$, $b' = b_1b$ in lemma 6.3, we obtain $f(cb) = h(cb)$.
Finally, take $a' = aa_1$ and $b' = b_1b$ in lemma 6.3. We obtain
$f(acb) = h(acb)$. $\square $

\it Remark 6.1. \rm Let $a_i$ be simple loops in $\Sigma_{0,5}$ so that
$a_i \perp_0 a_{i+1}$ and $|a_i \cap a_j| = 0$ for $|i-j| \geq 2$
as shown in figure 24(c). Then $a_i a_{i+1} \cong a_{i+2} a_{i+3} a_{i+4}$.
This seems to be an important relation on $\Cal S$. Indeed, let
$F_0 =\{ a_i, a_i a_{i+1}, b_i: i=1,..., 5\}$ where $b_i$'s are
the boundary components. Then the proof of lemma 6.1 shows that
$\Cal S(\Sigma_{0,5}) = \cup_{n=0}^{\infty} F_n$ where
$F_{n+1} = F_n \cup \{ \alpha | \alpha = \beta \gamma$ with
$\beta \perp_0 \gamma$, $\beta, \gamma, \gamma \beta$ and
the four components of $\partial N(\beta \cup \gamma)$ are in $F_n$\}.
The corresponding  curves for $\Sigma_{1,2}$ are as shown in figure 24.

\S7. Proofs of Theorem 1 and the Corollary 

To prove theorem 1 for  $\Sigma$ = $\Sigma_{g,r}$ with $|\Sigma| \geq 6$,
we decompose  $\Sigma$ = $X \cup Y$ as in \S5 so that $ X \cong \Sigma_{1,1}$
or $\Sigma_{0,4}$ and $\partial X \cap int(\Sigma)$ is a separating
simple loop $d$. By the reduction process in \S5, we construct a
measured lamination $m \in \Cal ML(\Sigma)$ so that $f = I_m$ on
the subset $\Cal S(X) \cup \Cal S(Y)$. To show that $f =I_m$,
by proposition 5.1 for $\delta = [d]$, it suffices to show that $f(\alpha)
=I_m(\alpha)$ for $\alpha \perp_0 [d]$. Take $x \in \alpha$ so that
$|x \cap d | =2$ and consider the incompressible surface
$\Sigma' = X \cup N(x)$. Then $\Sigma' \cong \Sigma_{0,5}$ or $\Sigma_{1,2}$.
Let $Y' = \Sigma' \cap Y$. Then $\Sigma' = X \cup Y'$ with $X \cap Y'
= X \cap Y \cong \Sigma_{0,3}$.  In particular $Y'$ is incompressible in
$\Sigma'$. Now consider $f|_{\Cal S(\Sigma')}$ and $I_m|_{\Cal S(\Sigma')}$.
They have the same values at elements in $\Cal S(X) \cup \Cal S(Y')$. Thus
by theorem 1  for $\Sigma'$ and lemma 2.1, we have  $f|_{\Cal S(\Sigma')} =$
$I_m|_{\Cal S(\Sigma')}$. In particular $f(\alpha) = I_m(\alpha)$.

To prove the corollary in \S1 for surface $\Sigma_{g,r}$ with
$\partial \Sigma = b_1 \cup ... \cup b_r$, we choose a Fenchel-Nielsen system
$\alpha = \alpha_1.... \alpha_n $ for $\Sigma$ where $n = 3g+r-3$. For each
index $i$, choose $\beta_i \in \Cal S'(\Sigma)$ so that
$I(\beta_i, \alpha_j) = 0$ for $j \neq i$ and $\beta_i \perp \alpha_i$
or $\beta_i \perp_0 \alpha_i$. We call the set $F =\{\alpha_i,
\beta_i, \alpha_i \beta_i, b_j: i=1,...,n, j=1,...r\}$ 
a \it Thurston basis \rm of the
measured lamination space. It is shown in [FLP] that the map $\tau_F
: \Cal ML(\Sigma) \to R_{\geq 0}^F$ sending $m$ to $I_m|_F$ is an
embedding (In [FPL], the set $F$ is taken to be
\{$\alpha_i, \beta_i, \alpha_i \beta_i^2, b_j$\}. But the proof
works  for our case as well). 
We shall show that the image of $\tau_F$ is a  semi-real
algebraic polyhedron by induction on $|\Sigma|$. By theorem 2, the
result holds for $|\Sigma| =4$. 
Now if $|\Sigma| \geq 5$,
we decompose $\Sigma = X \cup Y$ so that (1) $3 < |X|, |Y| < |\Sigma|$,
(2) $X \cap Y \cong \Sigma_{0,3}$ and (3) the components $a_1, a_2,
a_3$ of $\partial (X \cap Y)$ represent elements, say, $\alpha_1$, $\alpha_2$,
$\alpha_3$  in $F$ ($\alpha_2$ may be the same as $\alpha_3$).
Let $F_X = F \cap \Cal S(X)$ and $F_Y = F \cap \Cal S(Y)$.  There are
two possibilities: either $\alpha_1$, $\alpha_2$, $\alpha_3$ are
pairwise distinct or $\alpha_2 = \alpha_3$ ($\neq \alpha_1$). In the
first case, then
$F_X$ and $F_Y$ are Thurston bases for $X$ and $Y$ by condition (3) and
the definition. Let $\tau_{F_X}(m) = (x_1, ..., x_k)$ and $\tau_{F_Y}(m)
=(y_1,..., y_l)$ so that $x_i = I_m(\alpha_i)$, $y_i = I_m(\alpha_i)$
for $i=1,2,3.$ By the induction hypothesis, both images $Imag(\tau_{F_X})$
and $Imag(\tau_{F_Y})$ are semi-real algebraic polyhedrons. Now by
lemma 2.1, each $m \in \Cal ML(\Sigma)$ is determined by its restriction
on $X$ and $Y$.  Thus $Imag(\tau_F) =\{(x_1,..., x_k; y_1,..., y_l)
\in Imag(\tau_{F_X}) \times Imag(\tau_{F_Y})$: $x_i = y_i$, $i=1,2,3$\}.
Thus the result follows by the induction hypothesis.
In the second case that $\alpha_2 =\alpha_3$, one of the surfaces $X, Y$,
say, $X$ is $\Sigma_{1,1}$. Then $F_X$ is a Thurston basis for $X$ and
$F_Y \cup \{[a_2], [a_3]\}$ is a Thurston basis for $Y$. 
Let $\tau_{F_X}(m) = (x_1, ..., x_k)$ and $\tau_{F_Y}(m)
=(y_1,..., y_l)$ so that $x_i = I_m(\alpha_i)$, $i=1,2,$ and
$y_1 = I_m(\alpha_1),$ $y_2 = I_m([a_2])$, $y_3 = I_m([a_3])$. By the
same argument as above (using $x_1 = y_1$, $x_2 = y_2 = y_3$), 
the result follows.

\S8. Proofs of Results in Section 2 and Some Questions

8.1. The Proofs

{\bf Lemma  8.1.} \it (a) If $a$ and $b$ are curve systems with $|a \cap b|
= I(a,b)$, then the disjoint union $ab$ is a curve system.

(b) Suppose $a$, $b$ and $c$ are curve systems in $\Sigma$ 
so that $|a \cap b| = I(a,b)$, $|b \cap c| = I(b,c)$, $|c \cap a|
=I(c,a)$ and $| a \cap b \cap c| = 0$. If there is no contractible region in
$\Sigma -(a \cup b \cup c)$ which is either bounded by three arcs
in $a$, $b$ and $c$ respectively, or by four arcs in $a$, $b$, $c$ and
$\partial \Sigma$ respectively (see figure 29(a)), then $|c \cap ab| =
I(c, ab)$.

\rm

{\bf Proof.} (a) If $ab$ is not a curve system, then there exists either
(1) a simple closed curve $s$ in $ab$ and an annulus $D$ with $\partial
D$ = $s \cup d$ where $d$ is a boundary component of $\Sigma$ or
(2) a simple closed curve or a proper arc $s$ in $ab$ and a disc $D$
in $\Sigma$ so that either (2.1) $\partial D = s$ or (2.2) $\partial D$
$= s \cup d$ where $ s \cap d = \partial s = \partial d$ and $d$ is an
arc in $\partial \Sigma$. By replacing $s$ and finding another component
of $ab$ in $int(D)$ if necessary, we may assume that $ab \cap int (D)$
$= \emptyset$. Take a small regular neighborhood $N(a \cup b)$ of $a \cup b$ to be
$N(a) \cup N(b)$. We assume the resolutions are taken place inside
$N(a) \cap N(b)$. Thus $int(D)$ contains a finite number of
connected components $R_0, R_1,..., R_n$ of $\Sigma -$int$(N(a) \cup N(b)$),
where $R_i \neq \emptyset$, and $R_i \cap \partial \Sigma = \emptyset$,
 for $i \geq 1$, and $R_0 = \emptyset$ in case $D$ is a disc in $ int(\Sigma)$,
and  $R_0$ is the region which intersects $\partial \Sigma$ in the other cases.
Furthermore,  $R_0$ is a disc if $D$ is a disc intersecting $\partial \Sigma$ and
is an annulus if $D$ is an annulus. Each region $R_i$ ($i \geq 1$) is a disc
since otherwise there would be at least two boundary components of $R_i$ in
$int(D)$. 
This would contradict the assumption that $ int(D) \cap ab = 
\emptyset$. Call a point in $
\partial N(a) \cap \partial N(b)$ a \it corner \rm of $N(a \cup b)$. Each
point $p$ in $a \cap b$ corresponds to four corners in $\partial N(p)$
where $N(p)$ is the connected component of $N(a) \cap N(b)$
containing $p$. Join opposite corners in $\partial N(p)$ by
an arc in $int(N(p))$ so that it avoids one of the resolutions of $a \cap b$
at $p$. We call the arc a \it bridge \rm between the corners.
A corner of $\partial N(a \cup b)$ in a region $R_i$
is called a \it vertex \rm of $R_i$. Vertices of $R_i$ decompose $\partial R_i$
into \it edges. \rm Each edge is either in $\partial N(a)$, or in $\partial N(b)$, or
in $\partial \Sigma$. There is at most one edge which
is in $R_0$.  If two edges have a vertex in common, they cannot be both
in $N(a)$ (resp. in $N(b)$). Thus for $i \geq 1$, there are even number
of edges in $R_i$. Each region $R_i$ with $i \geq 1$, must have at least
four edges since $|a \cap b| = I(a,b)$ (if there were regions with only
two edges, then the region provides a Whitney disc for $a \cup b$).
More importantly, the definition of the resolution implies the
following alternating  principle: if $v$ and $v'$ are two 
vertices joint by an edge in $R_i$ so that the edge is either in $N(a)$
or in $N(b)$ then exactly
one of the bridges from $v$ or $v'$ still lies in $D$ (see figure 28(b)).

\midspace{0.1cm}
\centerline{\epsfbox{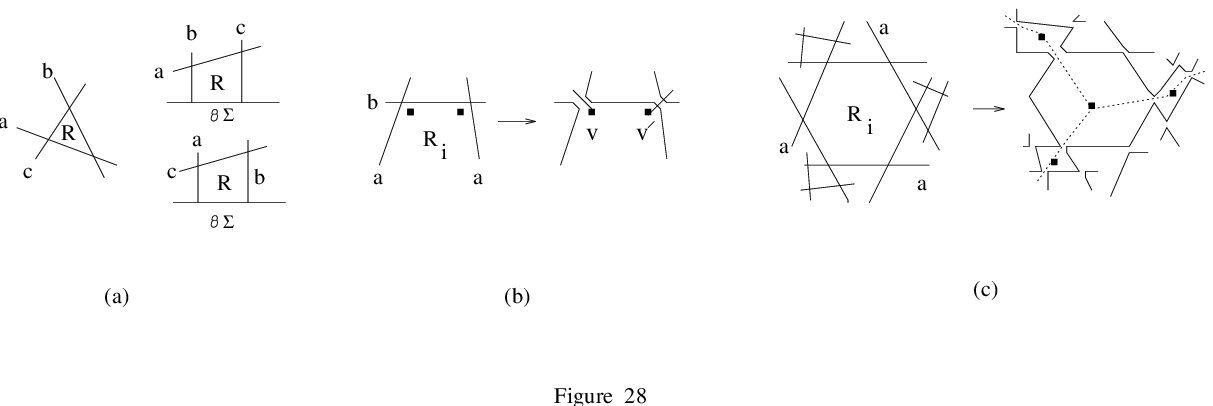}}
\midspace{0.1cm}

Form a graph $G$ in $D$ by putting a 0-cell in each
$int(R_i)$. Joint two 0-cells of $int(R_i)$ and $int(R_j)$ by a
1-cell in $D$ if there are opposite vertices in $R_i$ and $R_j$
so that their bridge is in $D$ (the 1-cell is an extension of the bridge).
These 1-cells are chosen to be pairwise disjoint except at the end points.
By the construction, if $D$ is a disc, the graph $G$ is homotopic to $D$
since each region $R_i$ is a disc; if $D$ is an annulus, the region $R_0$
is an annulus, thus the graph $G$ is  again homotopic to a disc. In both
cases, $G$ is a tree. Therefore either $G$ is a point or $G$ contains
two 0-cells of valency one. However by the construction, each region $R_i$
($i \geq 1$) has at least four edges and thus corresponds to a
 0-cell of valency at
least two by the alternating principle. Thus the graph $G$ must be a point.
Therefore, there is only one region $R_0$ which has at most one vertex
by the alternating principle. This contradicts the condition that 
$|a \cap b| = I(a, b)$.

(b)  Suppose the result is false. Then there is a disc $D \subset \Sigma$ so 
that either (1) $\partial D$ is a union of two arcs $s$ and $t$ with
$s \cap t = \partial s = \partial t$, $s \subset c$ and $t \subset ab$,
or (2) $\partial D$ is a union of three arcs $s, t, u$ so that each pair
of arcs intersect at one end point and $s \subset c$, $t \subset 
ab$, and $u \subset \partial \Sigma$. By taking the inner most disc if 
necessary, we may assume that $int(D) \cap (c \cup ab) = \emptyset $. 
Let $N(ab) = N(a) \cup N(b)$, $N(a \cap b) = N(a) \cap N(b)$, and $R_0, R_1,
..., R_n$ be the set of components of $\Sigma -( c \cup N(a) \cup N(b))$
which are contained in $D$. We set $R_0$ to be the region so that $R_0 \cap
c \neq \emptyset$. Then $R_0 \cap u \neq \emptyset$ if $u \neq \emptyset$.
Furthermore, $R_i \cap (c \cup \partial \Sigma) = \emptyset$ for $i \geq 1$.
By the assumption that int$(D) \cap (c \cup ab) = \emptyset $, each region $R_i$ is a disc.
Use the same argument as in (a), each region $R_i$ ($i \geq 1$) has at
least four sides and adjacent vertices in $\partial R_i$  ($i \geq 0$)
satisfy the alternating principle.  Form the same type of graph $G$ in $D$ based on the combinatorics of the regions $R_i$ as in (a). Since each region $R_i$
is contractible, the graph $G$ is a tree. Thus $G$ is either a point or
contains two vertices of valency one. The later case is impossible by
the alternating property. Thus $G$ is a point. Thus, there is only one region
$R_0$ in $D$ which has exactly one vertex. This is equivalent to the
condition that there is a contractible region in
$\Sigma -(a \cup b \cup c)$ which is bounded by three arcs in $a$, $b$, and
$c$, or by four arcs in $a$, $b$, $c$, and $\partial \Sigma$. Thus we obtain
a contradiction.
$\square$

{\bf Proposition 2.1.} \it The multiplication $\Cal CS(\Sigma)$$\times$$\Cal CS( \Sigma)$ $\to$ $\Cal CS(\Sigma)$ sends
$\Cal CS_0(\Sigma)$$\times$ $\Cal CS_0(\Sigma)$ to $\Cal CS_0(\Sigma)$ and 
satisfies the following properties.

(a) It is invariant under the action of the orientation preserving
homeomorphisms.

(b) If $I(\alpha, \beta)$ =0, then $\alpha \beta = \beta \alpha$. Conversely, if
 
$\alpha \beta = \beta \alpha$ and $\alpha \in$ $\Cal CS_0(\Sigma)$, then 
$I(\alpha, \beta)$=0. 

(c) If $\alpha$ $\in$ $\Cal CS_0(\Sigma)$, $\beta$ $\in$ $\Cal CS(\Sigma)$, then
 $I(\alpha, \alpha \beta)=I
(\alpha, \beta \alpha) = I(\alpha, \beta)$ and $\alpha(\beta \alpha)=
(\alpha \beta) \alpha$. If in addition that each component of
$\alpha$ intersects $\beta$, then $\alpha(\beta \alpha) = \beta$.

(d) If $[c_i] \in$ $\Cal CS(\Sigma)$ so that $|c_i \cap c_j| = I(c_i, c_j)$ for
$i,j=1,2,3$, $i \neq j$, $|c_1 \cap c_2 \cap c_3| =0$,
and there is no contractible region in $\Sigma -(c_1 \cup c_2 \cup c
_3)$
bounded by three arcs in $c_1$, $c_2$, $c_3$, then $[c_1]([c_2][c_3])=
([c_1][c_2])[c_3]$.

(e) For any positive integer $k$, $(\alpha ^k \beta^k) =(\alpha \beta)^k$.

(f) If $\alpha$ is the isotopy class of a simple closed curve, then the
positive Dehn twist along $\alpha$ sends $\beta$ to $\alpha^k \beta$ 
where $k =$ $I(\alpha, \beta)$.
\rm 

{\bf Proof.} Properties (a), (e) and (f) follow from the definition 
(see figure 29(a)).
Property (d) follows from lemma 8.1(b). Indeed, by the lemma, 
both $([c_1][c_2])
[c_3]$ and $[c_1]([c_2][c_3])$ are obtained by simultaneously resolving
all intersection points in $c_1 \cup c_2 \cup c_3$ from $c_1$ to $c_2$,
$c_2$ to $c_3$, and $c_1$ to $c_3$. To see (c), take $a$  and $a'$ to be
in $\alpha$ with $a \cap a' = \emptyset$ (two nearby parallel copies), and
$b \in \beta$ with $|a \cap b | = |a' \cap b| = I(a,b)$. Then, since
$\alpha $ is closed, $a, a', $ and $b$ satisfy the condition in lemma 8.1(b). 
Thus  $\alpha (\beta \alpha) = (\alpha \beta)\alpha$ follows. Also
by lemma 8.1(b), $I(\alpha, \alpha \beta) =  |a \cap a'b| = |a \cap b| =
I(\alpha, \beta)$ where $|a \cap a'b| = |a \cap b|$ follows from
the definition. 
The equality  $I(\alpha, \beta \alpha) = I(\alpha, \beta)$
follows similarly. 
If each component of $\alpha$ intersects $\beta$, then figure 29(b) shows
that $\alpha (\beta \alpha) = \beta$. Indeed, it suffices to
consider two adjacent intersection points $P_1$, $P_2$ along a 
component of $a$. Figure 29(b) shows that the multiplication
$a(b a')$ is the  same as finger moves on $b$. Thus 
 $\alpha(\beta \alpha) = \beta$. It remains to show (b). Clearly
if $I(\alpha, \beta) = 0$, then  $\alpha \beta$ = $\beta \alpha$.
Conversely, suppose $\alpha \in CS_0(\Sigma)$ and $\beta \in CS(\Sigma)$
with  $\alpha \beta$ = $\beta \alpha$. We decompose $\alpha$ as a
disjoint union $\alpha_1 \alpha_2$ where $I(\alpha_1, \beta) =0$ and
each component of $\alpha_2$ intersects $\beta$. Now since
$\alpha_1$ is disjoint from both $\alpha_2$ and $\beta$, we have $\beta
(\alpha_1 \alpha_2) = \alpha_1 (\beta \alpha_2)$.
Thus, by  $\alpha \beta$ = $\beta \alpha$, we obtain $\alpha_2 \beta 
= \beta \alpha_2$. Since each component of $\alpha_2$ intersects $\beta$, by 
property (c), $\beta = \alpha_2(\beta \alpha_2) = \alpha_2(\alpha_2 \beta)
= (\alpha_2)^2 \beta$ where the last equality follows from  property (d).
Now by property (c), $I(\beta, \beta) =
I(\beta, (\alpha_2^2) \beta)$  $= I(\beta, \alpha_2^2)$ = $2 I(\beta, \alpha_2)$$\neq 0$. This is a contradiction.
$\square$

\midspace{0.1cm}
\centerline{\epsfbox{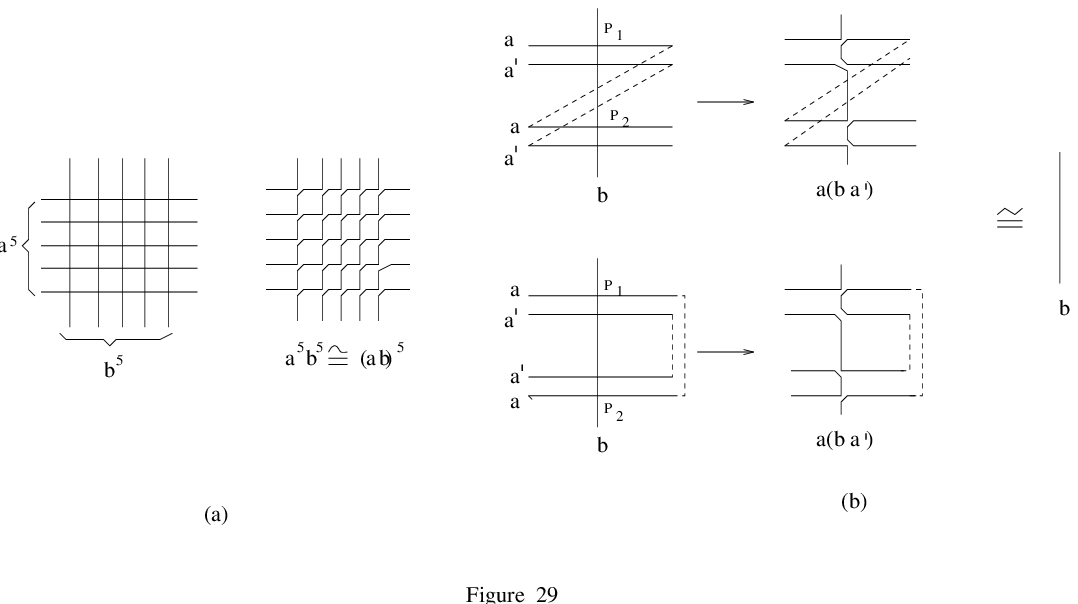}}
\midspace{0.1cm}

\it Remark 8.1. \rm   Properties (b), (c) and (d) are similar to the
commutative, the inverse, and the  associative laws in group theory.
Indeed, if  each component of  a curve system 
$\alpha \in \Cal CS_0$ intersects  both $\beta$ and $\gamma$
and $\beta \alpha = \gamma \alpha$, then (c) implies that $\beta
= \alpha(\beta \alpha) = \alpha (\gamma \alpha) = \gamma$.

8.2. Some observations and questions

We begin with a lemma.

{\bf Lemma 8.2.} \it Suppose
 $\alpha, \beta \in \Cal CS_0(\Sigma)$, then $I(\alpha \beta, \beta \alpha)
= 2I(\alpha, \beta)$. In particular, we have (a) $(\alpha \beta) (\beta \alpha) 
= \beta ^2 \gamma$ where $\gamma$ is disjoint from both $\alpha$
and $\beta$, and (b) $I(\delta, \alpha \beta) + I(\delta, \beta \alpha)
\geq 2 I(\delta, \beta)$.\rm

Proof.  
Choose $x, x' \in \alpha$, $y, y' \in \beta$ so that  $|x \cap y| = I(x,y)$
and $x'$, $y'$ are parallel copies of $x, y$. Then 
by the definition of the multiplication, there are no bi-gons in
$xy \cup y'x'$ and $| xy \cap y'x'| = 2I(x,y)|$. Thus
$ (xy)(y'x')$ is a representative of $(\alpha \beta)(\beta \alpha)$. An
easy calculation shows that $(xy)(y'x') \cong y^2 z$ where $z$ consists
of components of $x$ disjoint from $y$. Thus the lemma follows. $\square$

Call the
set \{$\alpha ^n \beta$: $n \in \bold Z$\} a \it horocycle \rm in the
space $\Cal CS_0(\Sigma)$.  It follows from the lemma 8.2 (b)
that $I_{\delta}$ is convex along horocycles. 
An element $\alpha \in \Cal CS_0(\Sigma)$ is called
\it maximal \rm if the only simple loops or arcs which 
 are disjoint from $\alpha$
are the components of $\alpha$.  
 By proposition 2.1(c),   for any two elements $\alpha$, $\beta$ 
in $\Cal CS_0(\Sigma)$, there is a horocycle  containing two maximal
 elements $\alpha'$ 
and $\beta'$ so that $I(\alpha, \alpha') = 0$ and $I(\beta, \beta') =0$.
Indeed, the horocycles  is of the form $\alpha' (\beta' \alpha')^n$.
The analogous  fact for measured laminations with respect
to the  extension of the earthquake was proved by Bonahon [Bo3]
and Papadopoulos [Pa1].

We close  with two questions. Call a function in
$n$ variables $x_1,..., x_n$ algebraically piecewise linear, or simply,
APL, if it is obtained from $x_i$'s by finite number of summations,
multiplications over $\bold Q$, and the absolute value $|.|$ operation.

Question 3. Does the  multiplication on $\Cal CS(\Sigma)$ extend 
to a  multiplication  on $\Cal ML(\Sigma)$ which is
APL with respect to the Thurston coordinate $\tau_F$? 

Question 4. Is  the intersection number function on  $\Cal ML(\Sigma)
\times \Cal ML(\Sigma) \to \bold R$ sending $(m, m')$ to  $I(m,m')$
a APL map with respect to the Thurston coordinates?

See [Bo2], [Bo3], [Pa1], [Pa2], [Re]  and the references cited therein
for more information  on the intersection numbers and
the earthquakes in the measured laminations space.

\centerline{\bf Reference}

[Bi] Birman, J.S.,  Braids, links, and mapping class groups.
Ann. of Math. Stud., 82, Princeton Univ. Press, Princeton, NJ, 1975

[Bo1] Bonahon, F.: The geometry  of Teichm\"uller space via geodesic
currents. Invent. Math.  92 (1988), 139-162.

[Bo2] Bonahon, F.: Bouts des vari\'et\'es hyperboliques de dimension $3$. 
 Ann. of Math. (2) 124 (1986), no. 1, 71--158.

[Bo3] Bonahon, F.: Earthquakes on Riemann surfaces and on 
measured  geod- \newline 
esic laminations. Trans. Amer. Math. Soc. 330 (1992), no. 1, 69--95. 

[Br] Brumfiel, G. W.: The real spectrum compactification of
Teichm\"uller space, Contemp. Math.,  74, AMS, (1988), 51-75.

[CM] Culler, M., Morgan, J.: Group actions on {\bf R}-trees. Proc. London Math. Soc.  55 (1987), no. 3, 571--604.

[CS] Culler, M., Shalen, P.: Varieties of group representations and
splittings of 3-manifolds. Ann. Math.  117 (1983), 109-146.

[De] Dehn, M.: Papers on group theory and topology. J. Stillwell (eds.).
 Springer-Verlag, Berlin-New York, 1987.

[FK] Fricke, R., Klein, F.: Vorlesungen  \"uber die Theorie der
Automorphen Functionen. Teubner, Leipizig, 1897-1912.

[FLP] Fathi, A., Laudenbach, F., Poenaru, V.: Travaux de Thurston sur les
surfaces. Ast\'erisque  66-67, Soci\'et\'e Math\'ematique de France, 1979.

[Ga] Gardiner, F.: Teichm\"uller theory and quadratic differentials,
John Wiley \& Sons, New York, 1987.

[GiM] Gilman, J., Maskit, B.: An algorithm for $2$-generator Fuchsian
groups. Michigan Math. J. 38 (1991), no. 1, 13-32.

[Go] Goldman, W.: Topological components of spaces of representations.
Invent. Math. 93 (1988) no.3, 557-607.

[GoM] Gonzlez-Acua, F., Montesinos-Amilibia, J.: On the
character variety of group representations in SL(2, $\bold C$) and
PSL(2, $\bold C$). Math. Z. 214 (1993), no. 4, 627--652.

[Har] Harer, J.: The second homology group of the mapping class group
of an orientable surface. Invent. Math. 72 (1083), 221-239.

[Hat] Hatcher, A.: Measured lamination spaces for surfaces. Topology Appl.
30 (1988), 63-88.

[Ho] Horowitz, R.: Characters of free groups represented in the 2-dimensional
special linear group. Comm. Pure and App. Math. 25 (1972), 635-649.

[HM] Hubbard, J.,  Masur, H.: Quadratic differentials and
foliations. Acta Math. 142 (1979), no. 3-4, 221--274.

[HT] Hatcher, A. Thurston, W.: A presentation for the mapping class group
of a closed orientable surface. Topology 19 (1980), 221-237.

[Ke] Keen, L.: Intrinsic moduli
on Riemann surfaces. Ann. Math.  84 (1966), 405-420

[Ker1] Kerckhoff, S.: The asymptotic geometry of Teichmüller
space. Topology 19 (1980), no. 1, 23--41. 

[Ker2] Kerckhoff, S.: The Nielsen realization problem. Ann. of
Math.  117 (1983), no. 2, 235--265.

[Li] Lickorish, R.: A representation of oriented 
combinatorial 3-manifolds. Ann.  Math.  72 (1962), 531-540

[Lu1] Luo, F.:  Geodesic length functions and Teichm\"uller spaces, 
J. Diff. Geom., to appear.

[Lu2] Luo, F.:  Geodesic length functions and Teichm\"uller spaces.
Electronic Research Announcement, AMS, 2 (1996),  no. 1, 34-41.

[Lu3] Luo, F.:  Automorphisms of the complex of curves, preprint.

[Mag] Magnus, W.: Rings of Fricke characters and automorphism groups of
free groups. Math. Zeit. 170 (1980), 91-103.

[Mo] Mosher, L.: Tiling the projective foliation space of a punctured surface.
Trans. Amer. Math. Soc. 306 (1988), 1-70.

[MS] Morgan, J., and Shalen, P.: Valuations, trees, and degenerations of
hyperbolic structures, I. Ann. of Math., 120 (1984), 401-476.

[Pa1] Papadopoulos, A.: On Thurston's boundary of Teichm\"uller
space and the extension of earthquakes. Topology Appl. 41 (1991), no. 3,
147--177.

[Pa2] Papadopoulos, A.: Geometric intersection functions and
Hamiltonian flows on the space of measured foliations on a surface. Pacific
J. Math. 124 (1986), no. 2, 375--402.

[Par] Parry, W.: Axioms for translation length functions.
Arboreal group theory (Berkeley, CA, 1988), 295--330, Math. Sci. Res. Inst.
Publ., 19, Springer, New York, 1991.

[PH] Penner, R.,  Harer, J.: Combinatorics of train tracks.
Annals of Mathematics Studies, 125. Princeton University Press, Princeton,
NJ, 1992.

[Re] Rees, M.: An alternative approach to the ergodic theory of
measured foliations on surfaces. Ergodic Theory Dynamical Systems 1 (1981),
no. 4, 461-488. 

[Sc] Schneps, L.: The Grothendieck theory of dessins d'enfants. 
 Cambridge University Press, Cambridge, New York, 1994.

[Th] Thurston, W.: On the geometry and dynamics of diffeomorphisms of
surfaces. Bul. Amer. Math. Soc. 19 (1988) no 2, 417-438.       

[Wa] Wajnryb, B.: A simple presentation for the mapping class group
of an orientable surface. Israel J. Math. 45 (1983), 157-174.

[Wo] Wolpert, S.: Geodesic length functions and the Nielsen problem.
J. Diff. Geom.  25 (1987), 275-296

\end